\documentclass[12pt]{article}
\usepackage{amsfonts,amssymb,amsmath}
\usepackage[T2A]{fontenc}  
\usepackage[cp1251]{inputenc}
\usepackage[english,russian]{babel}
\usepackage{geometry}
\newtheorem{theorem}{Theorem}
\newtheorem{lemma}[theorem]{Lemma}

\newcommand{\doc}{{\bf Proof. }}

\def\ds{\displaystyle}
\def\f{\frac}
\def\mb{\mbox}
\def\div{{\rm div}}
\def\C{\mathbb{C}}
\def\R{\mathbb{R}}

\def\K{\mathcal{K}}

\def\per{{\rm per}}
\def\nab{\nabla}
\def\n{\nabla}
\def\pa{\partial}

\def\ve{\varepsilon}
\def\e{\varepsilon}

\def\C0{C^\infty_0(\R^d)}

\def\ld{L^2(\R^d)}
\def\rd{\R^d}

\def\ild{\int_{\R^d}}

\def\iy{\int_{Y}}
\def\beq{\begin{equation}}
\def\eeq{\end{equation}}

\begin{document}
\title{Improved approximations of resolvents in
homogenization of fourth-order operators  with periodic coefficients 
}


\author{ S.\,E. Pastukhova}

\date{}
\maketitle


\begin{footnotesize} 
 In the whole space $R^d$, $d\ge 2$, we study homogenization of a divergence form elliptic fourth-order operator $A_\varepsilon$ with measurable
 $\varepsilon$-periodic coefficients, where $\varepsilon$ is a small parameter.  For the resolvent $(A_\varepsilon+1)^{-1}$, acting as an operator from $L^2$ to $H^2$, we find an approximation with remainder term of order $O(\varepsilon^2)$
 as $\varepsilon\to 0$. Relying on this result,
  we construct the resolvent approximation with remainder  of order $O(\e^3)$
 in the operator $L^2$-norm. 
 We employ  
 two-scale expansions that involve smoothing.

\end{footnotesize}
\bigskip

\section{Introduction}
\setcounter{theorem}{0} \setcounter{equation}{0}

  In the whole space $\rd$, $d\ge 2$, we consider a fourth-order differential operator (DO)
 \begin{equation}\label{0}
A_\e=\sum\limits_{i,j,s,t}D_{ij}
(a_{ijst}(x/\e)
D_{st}), \quad D_{ij}=\f{\partial^2}{\partial x_i\partial x_j},\quad \e\in (0,1),\end{equation}
with summation over all indices  $i,j,s,t$ from 1 to d. Here, the 1-periodic symmetric measurable real-valued fourth-order tensor
$a(x){=}\{a_{ijst}(x)\}$
is subject to conditions of positive definiteness and boundedness which ensure that the operator is elliptic.

The well-known result in  \cite{ZKOK}, where $G$-convergence and homogenization of families of 
DO much more general than (\ref{0}) were studied, means, in particular, that  the resolvent of the operator $A_\e$ 
is close to the resolvent of some 
limit, or homogenized, operator $\hat{A}$
in the sense of strong operator convergence in $\ld$:
 $(A_\e+1)^{-1}f\to (\hat{A}+1)^{-1}f$ in $\ld$ as $\e\to 0$ for any $f\in \ld$.
 The limit operator
   $\hat{A}$ is of the same ellipticity class, as $A_\e$, but is essentially simpler:
 \begin{equation}\label{03}
 \hat{A}=\sum\limits_{i,j,s,t}D_{ij}(\hat{a}_{ijst}
D_{st}),
  \end{equation}
where the tensor $\hat{a}=\{\hat{a}_{ijst}\}$ is constant and may be calculated in terms of solutions of auxiliary
periodic problems. 

The stronger operator convergence of $A_\e$ to $\hat{A}$, namely, the uniform resolvent convergence was proved in \cite{V},
 \cite{AA16}, \cite{P16} and, moreover, the convergence    rate estimate with respect to the parameter  $\varepsilon$ was deduced:
 \begin{equation}\label{02}
\|(A_\e+1)^{-1}-(\hat{A}+1)^{-1}\|_{L^2(\rd)\to L^2(\rd)}\le C\e
 \end{equation}
 with the constant $C$ depending only on the spatial dimension and the ellipticity constant of the tensor  $a(x)$.
   In  \cite{AA16} and \cite{P16}, even the stronger convergence in the operator $(L^2(\rd)\to H^2(\rd))$-norm    was proved: 
\begin{equation}\label{04}
\|(A_\varepsilon+1)^{-1}-(\hat{A}+1)^{-1}-\varepsilon^2 \mathcal{K}_\varepsilon\|_{L^2 (\mathbb{R}^d)\to H^2 (\mathbb{R}^d)}\le C\varepsilon
\end{equation}
with  the constant $C$ of the same type as in (\ref{02}). Here, we have the approximation of the resolvent $(A_\varepsilon+1)^{-1}$ in the energy operator norm and this is the sum of the 
zeroth term $(\hat{A}+1)^{-1}$
and the correcting term $\varepsilon^2 \mathcal{K}_\varepsilon$.
 The operator $\mathcal{K}_\varepsilon$
is determined by solutions of auxiliary
periodic problems on the periodicity cell which are introduced to define the limit (or homogenized) tensor  $\hat{a}$.
Note that
\begin{equation}\label{05}
\|\varepsilon^2 \mathcal{K}_\varepsilon\|_{L^2 (\mathbb{R}^d)\to H^2 (\mathbb{R}^d)}{\le} C, \quad 
\| \mathcal{K}_\varepsilon\|_{L^2 (\mathbb{R}^d)\to L^2 (\mathbb{R}^d)}{\le} C,
\end{equation}
where the constants $C$ are of the same type, as in  (\ref{04}).
Thus, due to (\ref{05})$_2$, the estimate (\ref{02}) in the $L^2(\rd)$-operator norm  follows from 
(\ref{04}) if the weaker operator norm is taken and the term $\varepsilon^2 \mathcal{K}_\varepsilon$ is omitted being nonessential with respect to the weaker norm.
By a more subtle use of an 
 $H^2$-estimate of   type (\ref{04}), it was showed in \cite{PMSb} that the order 
of accuracy in the  $L^2$-estimate (\ref{02}) can be improved 
up to $\e^2$: 
\begin{equation}\label{06}
\|(A_\varepsilon+1)^{-1}{-}(\hat{A}+1)^{-1}\|_{L^2 (\mathbb{R}^d)\to L^2 (\mathbb{R}^d)}\le C\varepsilon^2
\end{equation}
where the constant $C$ depends only on the spatial dimension and the ellipticity constant of the tensor  $a(x)$.
Assuming no symmetry conditions on the tensor  $a(x)$, the approximation of the resolvent $(A_\varepsilon+1)^{-1}$ in the $L^2(\rd)$-operator norm becomes more complicated (see \cite{P20}), namely,
\begin{equation}\label{07}
(A_\e+1)^{-1}=(\hat{A}+1)^{-1}+\e\K_1+O(\varepsilon^2), 
\end{equation}
where the operator $\mathcal{K}_1$ is independent of $\e$ in contrast with its counterpart in (\ref{04}).

In \cite{P21}, the approximation of the resolvent $(A_\varepsilon+1)^{-1}$ in the $(L^2(\rd){\to} H^2(\rd))$-operator norm was obtained with 
remainder term of order $O(\e^2)$. This approximation should be constructed as a sum
\begin{equation}\label{08}
(\hat{A}_\e+1)^{-1}+\e^2\K_2(\varepsilon)+\e^3\K_3(\varepsilon)
\end{equation}
of the zero term of a new type (in comparison with (\ref{04}), (\ref{06}) and (\ref{07})) and the two correctors.
Here, the operator
$\hat{A}_\e=\hat{A}+\e B$ is a singular perturbation of $\hat{A}$, where  $B$ is a fifth-order DO.
The operators  $\K_2(\e)$, $\K_3(\e)$ and $B$ are determined with the help of solutions to auxiliary problems on the cell of periodicity (see \S3),  the number of which increases compared to the case  of order $\e$ approximation in the same operator norm.
   The approximation (\ref{08}) is specified in (\ref{2.13})--(\ref{2.16}). We note here only that
   \begin{equation}\label{09}
\|\varepsilon^2 \mathcal{K}_2(\e)\|_{L^2 (\mathbb{R}^d)\to H^2 (\mathbb{R}^d)}{\le} C, \quad 
\|\varepsilon^2 \mathcal{K}_3(\e)\|_{L^2 (\mathbb{R}^d)\to H^2 (\mathbb{R}^d)}{\le} C,
\end{equation}
where the constant $C$ depends only on the spatial dimension and the ellipticity constant of the tensor  $a(x)$.

Possessing   
the approximation of the resolvent $(A_\varepsilon+1)^{-1}$ in the $(L^2(\rd){\to} H^2(\rd))$-operator norm with remainder term of order $\e^2$, we  construct the approximation of the resolvent $(A_\varepsilon+1)^{-1}$ in the $(L^2(\rd){\to} L^2(\rd))$-operator norm 
of order $\e^3$, and this is is the main goal of the paper. To this end, we use the  approach proposed in \cite{PMSb} and \cite{P20} to obtain the improved $L^2$-approximations of the resolvent of order $\e^2$ basing on its approximation  of order $\e$ in the energy operator norm. The exact formulation of the main result is given in Theorem \ref{th5.1}.


\section{$\!\!\!\!\!\!$. The statement of the problem}\label{sec2}
 \setcounter{theorem}{0} \setcounter{equation}{0}
 
{\bf 2.1.}
 We study the fourth-order equation in 
 $\rd$
 with parameter $\ve{\in}( 0,1)$:
\begin{equation}\label{2.1}
  \begin{array}{cc}
  u^\e\in  H^2(\rd),  \quad     A_\e u^\e+u^\e=
     f,  \quad f\in L^2(\rd),& \\
A_\e {=}D^* \,a_\e(x)D,\quad a_\e(x){=}a(x/\e). & \\
  \end{array}
\end{equation}
Here
 $a(y){=}\{a_{ijst}(y)\}$ is a measurable real-valued fourth-order tensor, periodic
with periodicity cell $Y{=}[-\f12,\f12)^d$ and satisfying the symmetry and ellipticity conditions 
\beq\label{2.2}
\ds{a_{ijst}=a_{stij},\quad a_{ijst}=a_{jist},}
\atop\ds{
\exists\lambda>0:\quad\xi\cdot\xi\le a(\cdot)\xi\cdot\xi\le\lambda^{-1}\xi\cdot\xi
}
\eeq
 for any symmetric matrix $\xi{=}\{\xi_{ij}\}$, where
$\xi\cdot\eta$ denotes the inner product of the matrices матриц $\xi{=}
\{\xi_{ij}\}$ and 
$\eta{=}\{\eta_{ij}\}$, i.e. 
$\xi\cdot\eta{=}\xi_{ij}\eta_{ij}$. 
Throughout the paper, summation from $1$
to $d$ over repeated indices
is assumed  (unless otherwise specified). The constant $\lambda$ is called the ellipticity constant of the tensor $a(x)$.

By the symmetry property in (\ref{2.2}),  the action of the tensor  $a$
on a symmetric matrix $\xi$
is the symmetric matrix $\eta{=}a\xi$,
and for arbitrary   matrix $\xi$ we have
$a\xi{=}a\xi^s$,
where $\xi^s{=}1/2(\xi{+}\xi^T)$ is the symmetric part of $\xi$, or $a(\xi{+}\xi^T)=2a\xi$.

We need to define DO in (\ref{2.1}). Let
$$
D\varphi=\n^2\varphi=
\left\lbrace D_{ij}\right\rbrace_{i,j=1}^{d}, \quad D_{ij}=\f{\partial^2\varphi}{\partial x_i\partial x_j},
$$  
so that the operator $D$ sends scalar functions into matrix functions. Then the adjoint operator
$D^*{=}\n^*\n^*{=}\div\,\div$ acts formally on a matrix function
 $\eta{=}\{\eta_{ij}\}$
by   
$$D^*\eta=D_{ij}\eta_{ij}=\f{\partial^2}{\partial x_i\partial x_j}\eta_{ij},$$  sending a matrix  into a scalar function.
The equality
$
D^*\eta{=}f, 
$
where $\eta{\in} L^2(\rd)^{d\times d}$ and $ f {\in} L^2(\rd)$,
is understood in the sense of distributions on $\rd$.
So, the expanded form of the operator $A_\e$
 is given by
 \[
 A_\e=D_{ij}\left( a_{ijst}(x/\e)D_{st}\right),\quad D_{ij}=\f{\partial^2}{\partial x_i\partial x_j},
 \]
which is the same as in  (\ref{0}), by summation convention;
a solution of equation 
 (\ref{2.1}) is understood in the sense of the integral identity
\[
 \ild\left( a_\e(x) D u^\e \cdot D\varphi+ u^\e \varphi\right)\, dx=
 \ild f \varphi  dx\quad \forall \varphi\in   C_0^\infty(\rd),
\]
where the test functions can be taken from  
$ H^2(\rd)$, by the closure property.
By the Riesz Representation Theorem for a linear functional on a Hilbert space, the solution to
 (\ref{2.1}) exists and is unique, and satisfies an estimate
\[
 \| u^\e\|_{H^2(\rd)}\le c(\lambda)
 \|f \|_{L^2(\rd)}.
\] 
It worth noting that we can define a norm on the space 
  $ H^2(\rd)$, equivalent to the standard one, by the relation
 \[
 \|\varphi\|_{H^2(\rd)}=\left(\int_{\rd}(|\varphi|^2+|D\varphi|^2)\,dx\right)^{1/2}.
 \]

\bigskip 
{\bf 2.2.}
We introduce the equation
\begin{equation}\label{2.6}
  \begin{array}{cc}
  u\in  H^2(\rd),  \quad     \hat{A} u+u=
     f,  \quad f\in L^2(\rd),& \\
\hat{A} {=}D^* \,\hat{a}\,D, & \\
  \end{array}
\end{equation}
with the constant fourth-order tensor $\hat{a}$ defined in
 (\ref{3.5}). The expanded form of $\hat{A}$
 is given by
 $
 \hat{A}=D_{ij} \hat{a}_{ijst}D_{st}, 
 $
which is the same as in    (\ref{03}). Since the coefficients in (\ref{2.6}) are constant and satisfy the ellipticity condition       of the type (\ref{2.2}), it is easy to obtain
the following elliptic estimate  for the solution of this problem:
\beq\label{2.5} 
\|u\|_{H^4(\rd)}\le C \|f\|_{\ld},\quad C=const(d,\lambda).
\eeq

According to 
\cite{ZKOK},   the solutions to (\ref{2.1}) and (\ref{2.6})
are conected with the convergence
\beq\label{2.7}
\lim_{\e\to 0}\|u^\e-u\|_{\ld}=0\quad \forall f\in L^2(\rd),
\eeq
which allows to call  (\ref{2.6}) the limit, or homogenized, equation  for (\ref{2.1}). Since
$u^\e=(A_\e+1)^{-1}f$ and $u=(\hat{A}+1)^{-1}f$, and the resolvents $(A_\e+1)^{-1}$ and $(\hat{A}+1)^{-1}$ are operators in в $\ld$, the relation (\ref{2.7})  means the strong resolvent convergence in $\ld$  of
$A_\e$ to the operator $\hat{A}$. The recent results on operator convergence of the resolvent $(A_\e+1)^{-1}$ (moreover, with  convergence rate estimates that may contain appropriate correctors), formulated in (\ref{02}), (\ref{06}), 
 (\ref{04}) and (\ref{07}), 
also allows us to consider   $\hat{A}$ as the limit operator.

Approximating the resolvent  
$(A_\e+1)^{-1}$  
in the operator $(L^2(\rd)\to H^2(\rd))$-norm, the same as in (\ref{04}), but with an accuracy of $\e^2$ order,
 it is reasonable to take, instead of $(\hat{A}+1)^{-1}$,
the resolvent $(\hat{A}_\e+1)^{-1}$ of the operator $\hat{A}_\e$, more complicated than 
  $\hat{A}$. Namely,
\beq\label{2.8}
\hat{A}_\e=D_{ij} \left(\hat{a}_{ijst}D_{st}+\e b_{ij}^{rst}D_{rst}\right), \quad D_{rst}=\f{\partial^3}{\partial x_r\partial x_s \partial x_t},
\eeq
or more shortly written in terms of the operators  
$D$
and $D^*$,
\beq\label{2.9}
\hat{A}_\e u=D^* \left(\hat{a}D u+\e b^{rst}D_{rst}u\right), 
\eeq
where the matrices $b^{rst}=\{b_{ij}^{rst}\}_{i,j}$ are defined in (\ref{3.11}) for all indices $r,s,t$ from 1 to d.

 We introduce the new version of the homogenized equation with the operator (\ref{2.8})
\beq\label{2.10}
(\hat{A}_\e+1)\hat{u}^\e=f,\quad f\in \ld.
\eeq 
Since the coefficients of the equation are constant, its solution can be obtained with the help of 
the Fourier transform. Equation
(\ref{2.10}) is written in the Fourier images as
\beq\label{2.11}
\left(1+\Lambda(\xi)-i\e\Lambda_0(\xi)\right) F[\hat{u}^\e]=F[f],\quad i=\sqrt{-1},
\eeq
where $f(x)\to F[f](\xi)$ is the Fourier transform and
$\Lambda(\xi){=}\hat{a}_{pqst}\xi_p\xi_q\xi_s\xi_t,$ $\Lambda_0(\xi){=}b_{pq}^{rst}\xi_p\xi_q\xi_r\xi_s\xi_t.$
The coefficients $\hat{a}_{pqst}$ and $b_{pq}^{rst}$ are real; therefore,
 from  (\ref{2.11})  the inequality
\[
\ild \left(1+\Lambda(\xi)\right)^2 |F[\hat{u}^\e]|^2\,d\xi\le \ild|F[f]|^2\,d\xi,
\]
follows, thanks to ellipticity property of $\hat{a}$ and the  Plancherel identity. Hence, the counterpart of the elliptic estimate (\ref{2.5}) is derived, namely,
\beq\label{2.12}
\|\hat{u}^\e\|_{H^4(\rd)}\le C \|f\|_{\ld},
\eeq
where the constant $C$ depends only on the dimension d and the ellipticity constant 
$\lambda$
 from (\ref{2.2}). 


\bigskip 
{\bf 2.3.} We take the approximation for the solution of   (\ref{2.1}) in the form
\beq\label{2.13}
\tilde{u}^\e(x)=u^{,\e}(x)+\e^2\,U_2^\e(x)+\e^3\,U_3^\e(x)
\eeq
with
\beq\label{2.14}
U_2^\e(x)=N^{ij}(x/\e)D_iD_{j}u^{,\e}(x),\quad U_3^\e(x)=N^{ijk}(x/\e)D_iD_{j}D_ku^{,\e}(x),\quad D_i=\f{\pa}{\pa x_i},
\eeq
\beq\label{2.15}
u^{,\e}(x)=\Theta^\e\, \hat{u}^\e(x),\quad \Theta^\e=S^\e S^\e S^\e,
\eeq
where $S^\e$  is Steklov`s smoothing,
 $N^{ij}(y)$, $N^{ijk}(y)$,  for all indices $i,j,k$ from l to d,  and $ \hat{u}^\e(x)$ are solutions of the problems (\ref{3.1}), (\ref{3.9}),
  and (\ref{2.10}), respectively. We recall that for
   $\varphi\in L^1_{loc}(\rd)$ the Steklov smoothing operator is defined by 
\begin{equation}\label{2.16}
(S^\e\varphi)(x)=\iy \varphi(x-\e\omega)\,d\omega,\quad Y=[-1/2,1/2)^d.
\end{equation}

\begin{theorem}\label{th2.1}
For the difference of the solution to (\ref{2.1}) and the function given in 
(\ref{2.13})--(\ref{2.16}),
the following estimate holds:
\begin{equation} \label{2.17}
 \|u^\e-\tilde{u}^\e\|_{H^2(\rd)}\le C\e^2 \|f\|_{L^2(\rd)},
\end{equation}
where the constant   $C$  depends only on the dimension d and the ellipticity constant 
$\lambda$
 from (\ref{2.2}). 
\end{theorem}

By properties of smoothing,  the first term $u^{,\e}$ in (\ref{2.13})
can be replaced with $ \hat{u}^\e(x)$ preserving  
 the order of the majorant in  (\ref{2.17}); 
 in operator language,  this form of (\ref{2.17}) implies the approximation (\ref{08}) of the resolvent $(A_\varepsilon+1)^{-1}$ in the operator energy norm.
Theorem \ref{th2.1} is actually proved in \cite{P21}, where smoothing in the formulas
(\ref{2.13})--(\ref{2.16}) is performed by the double Steklov operator $S^\e S^\e$, which can be replaced  with an arbitrary smoothing operator  (\ref{m.11}) provided that its kernel is even and sufficiently regular. Smoothing 
$\Theta^\e{=}(S^\e)^3{=}S^\e(S^\e)^2$ is taken in (\ref{2.13})--(\ref{2.16}) because special features of the operators $S^\e$ and $(S^\e)^2$ given in Lemma  \ref{LemM6} are essential when we construct
 $\e^3$ order $L^2$-approximations relying on  $\e^3$ order $H^2$-approximations.

\bigskip 
{\bf 2.4.} Our main purpose is to present 
 such a function $w^\e(x)$ that is be constructed in terms of the solutions of the homogenized equation and several periodicity cell problems, and satisfies the estimate
\begin{equation} \label{2.20}
 \|u^\e-w^\e\|_{L^2(\rd)}\le C\e^3 \|f\|_{L^2(\rd)},
\end{equation}
where  the constant $C$  depends only on the spatial dimension $d$ and the ellipticity constant $\lambda$ from (\ref{2.2}). In operator language, this means that we obtain the approximation of the resolvent $(A_\e+1)^{-1}$, which is sharper than that of  (\ref{06}) and (\ref{07}), namely, its sharpness is  of $\e^3$ order.
The approximation of the resolvent has 
the form
\begin{equation}\label{2.21}
(\hat{A}_\e+1)^{-1}+\e^2\K(\varepsilon),
\end{equation}
where the zero term $(\hat{A}_\e{+}1)^{-1}$  of a new type  (compared with (\ref{06}) and (\ref{07})) is supplemented by rather complicated corrector,  whose 
four-part structure 
is determined in Theorem \ref{th5.1}. Here,
$\hat{A}_\e=\hat{A}+\e B$  is the same DO with constant coefficients as in (\ref{08}) and (\ref{2.9}).
The following estimate holds:
\begin{equation}\label{2.22}
\|(\hat{A}_\e+1)^{-1}+\e^2\K(\varepsilon)-(A_\e+1)^{-1}\|_{L^2(\rd)\to L^2(\rd)}\le C\e^3,
\end{equation}
where the constant $C$  depends only on the spatial dimension $d$ and the ellipticity constant $\lambda$ from (\ref{2.2}).

We now give an outline of further sections. Sections 3 and 4 are of an auxiliary and preparatory nature: they contain necessary  information about  cell problems  
and   smoothing operators respectively.  
 The approximation (\ref{2.21}) is directly constructed with its justification 
  in Section 5. Theorem \ref{th2.1} and Lemma  \ref{Lem5.1} are proved 
  in Section 6;
these are the basic preliminaries 
 to constructing the approximation  
  (\ref{2.21}).  Remarks connected  with 
  this  approximation 
   are incorporated 
   in Section 7.
 
\section{$\!\!\!\!\!\!$. Cell problems}\label{sec3}
 \setcounter{theorem}{0} \setcounter{equation}{0}
\textbf{3.1.} We introduce the following problem on the periodicity cell $Y$:
\begin{equation} \label{3.1}
 N^{ij}\in \tilde{H}^2_{per}(Y),\quad D^*_y[a(y)(D_yN^{ij}(y)+e^{ij})]=0,\,\, i,j=1,\ldots,d.
\end{equation}
Here  the tensor $a(y)$ is the same as in Section 2.1; the entries of the matrix
 $e^{ij}=\{e^{ij}_{sh}\}_{s,h}$  are $e^{ij}_{sh}=
\delta^i_s\delta^j_h$, where
$\delta^i_s$ is the Kronecker delta; 
  $\tilde{H}^2_{per}(Y){=}\{\varphi{\in}{H}^2_{per}(Y): \langle \varphi\rangle{=}0\}$ is the Sobolev space of  1-periodic functions with zero mean
 $$\langle \varphi\rangle=\iy \varphi\,dy.$$  
 The equality
\begin{equation} \label{3.2}
 D^*_yb=0
 \end{equation}
with $b\in L_{per}^2(Y)^{d\times d}$ means that the integral identity
\begin{equation} \label{3.3}
\langle    b\cdot D\varphi\rangle=0\quad \forall\varphi\in C^\infty_{per}(Y),
\end{equation}
holds, where the test functions can be taken from the space $\tilde{H}^2_{per}(Y)$ by closure; matrices satisfying (\ref{3.2})  are called \textit{solenoidal}.

In the theory of homogenization it is well known that  (\ref{3.2}) can be understood also in the sense of distributions on $\rd$; that is, as an alternative to (\ref{3.3}) we can consider the integral identity
\[
\ild b\cdot D\varphi\,dx=0\quad \forall
\varphi\in C_0^\infty(\rd).
\] 

By (\ref{3.2}),  solutions of the above problem (\ref{3.1}) satisfies the identity
\begin{equation} \label{3.3a}
\langle aD N^{ij}\cdot D\varphi\rangle=-\langle ae^{ij}\cdot D\varphi\rangle\quad \forall \varphi\in \tilde{H}^2_{per}(Y).
\end{equation}
 The existence of a unique solution $N^{ij}$ to equation  (\ref{3.1}) is proved using Riesz`s theorem.
It is appropriate to note here that the norm in
  $\tilde{H}^2_{per}(Y)$,
   introduced by
$\|\varphi\|_{\tilde{H}^2_{per}(Y)}=\langle|D\varphi|^2\rangle^{1/2},$
 is equivalent to
the standard norm thanks to the Poincar\'e inequality 
\[\langle|v|^2\rangle\le C_P\langle|\n v|^2\rangle\quad
\forall v \in H^1(Y),\,\langle v\rangle=0,
\]
which is valid  for any function $\varphi\in \tilde{H}^2_{per}(Y)$ and its derivatives $D_j \varphi$, $1\le j\le d$.

 From the integral identity (\ref{3.3a}),  it follows that
   \begin{equation} \label{3.4}
 \|N^{ij}\|^2_{\tilde{H}^2_{per}(Y)}\le C,
 \quad i,j=1,\ldots,d,
\end{equation}
with the constant $C$ depending on the dimension  d and the ellipticity constant $\lambda$.

 The fourth-order tensor $\hat{a}$ 
 from (\ref{2.6}) is defined by
\begin{equation} \label{3.5}
 \hat{a}e^{ij}=\langle    a(\cdot)(D_yN^{ij}(\cdot)+e^{ij})\rangle,
 \quad i,j=1,\ldots,d.
\end{equation}
 The tensor $\hat{a}$ inherits the properties  (\ref{2.2}) of symmetry and positive definiteness,
 which is proved in the same way as a similar fact from homogenization of second-order divergent elliptic equations (see,   
  for instance, \cite{ZKOK}, \cite{BP}, \cite{JKO}).

\textbf{3.2.} The following lemma generalizes the well known assertion about the representation of periodic solenoidal vectors with zero mean as the divergence of a skew-symmetric matrix
(see, \cite{JKO}, 
  for instance).
Its proof can be found in \cite{AA16}, \cite{P16} and \cite{UMN}.
 \begin{lemma}\label{lem1}
Let $ g{=}\{g_{st}\}_{s,t}\in L^2_{\per}(Y)^{d\times d}$ be a symmetric matrix such that
$\langle g\rangle=0$, $ D_y^*g{=}0$.
Then there exists a family of 1-periodic matrices  
$ G^{st}{=}\{G^{st}_{km}\}_{k,m}{\in} \tilde{H}^2_{\per}(Y)^{d\times d}$, $s,t=1,\ldots,d$,
satisfying conditions\\
$\quad{}\quad$ (i) $ G^{st}_{km}=G^{st}_{mk}$; \\ 
$\quad{}\quad$ (ii) $ G^{st}_{km}=-G^{km}_{st}$;\\ 
$\quad{}\quad$   (iii) $\| G^{st}\|_{H^2(Y)^{d\times d}}\le c(d)
 \|g \|_{L^2(Y)^{d\times d}}$;\\
$\quad{}\quad$ (iv) $g_{st}=D^* G^{st}$\\
for all indices $s,t,k,m$ from 1 to d.
\end{lemma}


 Lemma \ref{lem1} can be applied to the matrix function
\begin{equation} \label{3.6}
g^{ij}:=a(D N^{ij}+e^{ij})-\hat{a} e^{ij}, 
\end{equation}
since (\ref{3.1}) and (\ref{3.5}) guarantee that
\begin{equation} \label{s}
D_y^*g^{ij}{=}0,\quad \langle g^{ij}\rangle{=}0.
\end{equation}
Herefrom, the entries of the matrix $g^{ij}=\{g^{ij}_{st}\}_{s,t}$
can de represented in terms of the matrix potential  $ G^{ij,st}=\{G^{ij,st}_{km}\}_{k,m}\in \tilde{H}^2_{per}(Y)^{d\times d}$ in such a way that
\begin{equation} \label{3.8}
g^{ij}_{st}=D^* G^{ij,st}, \quad G^{ij,st}_{km}=-G^{ij,km}_{st}, \quad \| G^{ij,st}\|_{H^2(Y)^{d\times d}}\le c(d,\lambda)
\end{equation}
for all indices  from 1 to d.

\textbf{3.3.} Possessing the solutions $N^{ij}$
to the problem (\ref{3.1}) and matrix potentials $G^{ij,st}$ in representation (\ref{3.8}), we introduce another family of cell problems
\begin{equation} \label{3.9}
 N^{ijk}\in \tilde{H}^2_{per}(Y),\quad D^*\left( aDN^{ijk}+2a(\n N^{ij}\times e^k)+2\partial_m G^{ij,km}\right)=0,\,\, i,j,k=1,\ldots,d.
\end{equation}
Here, 
$e^1,\ldots,e^d$ is the canonical basis  in $\rd$; $\partial_m=\f{\pa}{\pa y_m}$;
   the matrix $\alpha{\times}\beta{=}\{\alpha_p\beta_q\}_{p,q}$ is constructed from two $d$-dimmmensional vectors $\alpha=\{\alpha_p\}_p$ and $\beta=\{\beta_q\}_q$, we set
  $\alpha{=}\n N^{ij}$ and $\beta{=}e^k$.

  For any triple of indices $i,j,k$  from 1 to d, the problem  (\ref{3.9}) is uniquely solvable; being related to it, the 
  matrix
  \begin{equation} \label{3.10}
 g^{ijk}:= aDN^{ijk}+2a(\n N^{ij}\times e^k)+2\partial_m G^{ij,km}-b^{ijk} 
\end{equation}
with
  \begin{equation} \label{3.11}
 b^{ijk}:=\langle aDN^{ijk}+2a(\n N^{ij}\times e^k)\rangle 
\end{equation}
satisfies the conditions of Lemma
 \ref{lem1} because  
\begin{equation} \label{ss}
D_y^*g^{ijk}{=}0,\quad \langle g^{ijk}\rangle{=}0.
\end{equation}
Hence, there is the representation for the entries of  $ g^{ijk}=\{g^{ijk}_{st}\}_{s,t}$
in terms of the matrix potential $ G^{ijk,st}=\{G^{ijk,st}_{pq}\}_{p,q}\in \tilde{H}^2_{per}(Y)^{d\times d}$ such that
\begin{equation} \label{3.13}
g^{ijk}_{st}=D^* G^{ijk,st}, \quad G^{ijk,st}_{pq}=-G^{ijk,pq}_{st}, \quad \| G^{ijk,st}\|_{H^2(Y)^{d\times d}}\le c(d,\lambda)
\end{equation}
for all indices from 1 to d.

 Since $D^* G^{ij,km}{=}g^{ij}_{km}$
  and the  matrix $g^{ij}$
  is expressed in terms of $N^{ij}$
  in view of (\ref{3.6}), we can rewrite equation (\ref{3.9}) for $N^{ijk}$
in such a way that  its right hand side function is determined only by the function $N^{ij}$ and its derivatives, the tensor $a$ and the vector $e^k$. For further calculations,
 we prefer the equation for $N^{ijk}$ in display (\ref{3.9}).

\textbf{3.4.} As a corollary of  \ref{lem1}, we have the next 
 \begin{lemma}\label{lem2}
Let $ g(y){=}\{g_{st}(y)\}_{s,t}$ be a matrix satisfying the conditions of Lemma
  \ref{lem1}. Let 
$ G^{st}(y){=}\{G^{st}_{km}(y)\}_{k,m}\in \tilde{H}^2_{\per}(Y)^{d\times d}$, where $s,t=1,\ldots,d$, be a corresponding family of matrix potentials 
with properties (i)--(iv) listed in Lemma \ref{lem1}. 
Then 
\begin{equation} \label{3.14}
\ds{g_{st}(x/\e)\Phi(x)}
\atop\ds{
=\e^2D^*\left(  G^{st}(x/\e)\Phi(x)\right)-
\e^2\,  G^{st}(x/\e)\cdot D\Phi(x)-2\e\,(\div_y  G^{st})(x/\e)\cdot\n\Phi(x)}
\end{equation}
for any
$\Phi \in \C0$  and all indices $s,t$ from 1 to d.
 Here, the matrix 
 \begin{equation} \label{3.15}
M(x)= \{D^*(  G^{st}(x/\e)\Phi(x))\}_{s,t}
\end{equation}
is solenoidal, that is, 
$
D^* M=0 \, (\mb{in the sense of distributions on} \, \rd).
$
\end{lemma}
\doc 
By assumptions, we have 
$$g_{st}(y)=D^* G^{st}(y),\quad g_{st}(x/\e)\Phi(x)=D^*\left(\e^2 G^{st}(x/\e)\right)\Phi(x);
$$
therefore, (\ref{3.14}) holds thanks to the rule for applying the operator $D^*$ 
to the product of a symmetric matrix  $B=\{B_{ij}\}_{i,j}$
 and a scalar $\Phi$:
 \begin{equation} \label{3.16}
D^*(\Phi\,B)=\Phi\,D^*B+B\cdot D\Phi+2(\div B)\cdot \n \Phi,
\end{equation}
where $ \div B=\{\partial B_{ij}/\partial x_j \}_i$ is a vector, and the inner product both of matrices and vectors is denoted similarly by dot.

The fact that
 the matrix $M(x)$ is solenoidal is ensured
by identity
\[
\ild \varphi\, D^* M\,dx=0 \quad \forall \varphi\in \C0.
\]
The latter follows from the chain of integral identities
\[
\ild \varphi\, D^* M\,dx=
\ild D\varphi\cdot M\,dx\stackrel{(\ref{3.15})}=\ild D_{st}\varphi(x) D^*(  G^{st}(x/\e)\Phi(x))\,dx
\]
\[
=\ild D_{st}\varphi(x)  D_{km}( G_{km}^{st}(x/\e)\Phi(x))\,dx
=\ild D_{km}D_{st}\varphi(x)   G_{km}^{st}(x/\e)\Phi(x)\,dx
\]
and the pointwise equality
$D_{km}D_{st}\varphi(x)   G_{km}^{st}(x/\e){=}0$, which is true  
 by the skew-symmetry property
$G_{km}^{st}=-G_{st}^{km}$ of the matrix potentials $ G^{st}(y){=}\{G^{st}_{km}(y)\}_{k,m}$.

\section{$\!\!\!\!\!\!$. Smoothing operators}\label{sec4}
 \setcounter{theorem}{0} \setcounter{equation}{0} 
 
For brevity, we denote the  norm and the inner product  in $\ld$ as
 \[
\|\cdot\|=\|\cdot \|_{\ld},\quad
(\cdot \,,\cdot\,)=(\cdot \,,\cdot\,)_{\ld},
\]
 making no distinguish between spaces of scalar and vector functions.

 \textbf{4.1. Steklov`s smoothing operator.} 
 We recall the definition of Steklov`s smoothing operator 
 \begin{equation}\label{m.1}
S^\e\varphi(x)=\iy\varphi(x-\e\omega)\,d\omega,\quad Y=[-1/2,1/2)^d,
\end{equation}
and give, firstly, 
the simplest and well known properties:
  \begin{equation}\label{m.2}
\|S^\e\varphi\|\le\|\varphi\|,
\end{equation}
 \begin{equation}\label{m.3}
\|S^\e\varphi-\varphi\|\le (\sqrt{d}/2)\e\|\nab\varphi\|\quad \forall\varphi\in H^1(\R^d).
\end{equation}
We mention also the evident property $S^\e(\nab \varphi)=\nab 
(S^\e\varphi)$, which is used systematically further. The following properties are pivotal for our method (see their proof in 
\cite{UMN}--\cite{ZhP06}). 
  \begin{lemma}\label{LemM1} If $\varphi{\in}L^2(\R^d)$, $b{\in}L^2_\per(Y)$,
 and $b_\e(x){=}b(\e^{-1}x)$,
 then $b_\e S^\e\varphi\in L^2(\R^d)$ and
 \begin{equation}\label{m.4}
\|b_\e S^\e\varphi\|\le\langle b^2\rangle^{1/2}\|\varphi\|.
\end{equation}
\end{lemma} 

 \begin{lemma}\label{LemM2} If $b{\in}L^2_\per(\Box)$, $\langle b\rangle{=}0$, $b_\e(x){=}b(\e^{-1}x)$,
  $\varphi{\in}L^2(\rd)$,
 and $\psi{\in}H^1(\rd)$,
 then
 \begin{equation}\label{m.5}
( b_\e S^\e\varphi,\psi) \le C\e
\langle b^2\rangle^{1/2}\|\varphi\|\,\|\nab\psi\|,\quad C=const(d).
\end{equation}
\end{lemma}

The above estimates 
can be improved under higher regularity conditions. For example,  
 \begin{equation}\label{m.6}
\|S^\e\varphi-\varphi\|\le C\e^2\|\nab^2\varphi\|\quad \forall\varphi\in H^2(\R^d),\quad C=const(d),
\end{equation}
which 
will be even more specified for
$\varphi{\in} H^4(\R^d)$ in Section 4.3.
By duality, from (\ref{m.6}) we obtain 
\begin{equation}\label{m.7}
\|S^\e\varphi-\varphi\|_{H^{-2}(\rd)}\le C\e^2\|\varphi\|_{L^2(\rd)}\quad\forall\varphi\in L^2(\R^d),\quad C=const(d).
\end{equation}

 The $L^2$-form in (\ref{m.5})  
 has a larger smallness order in the following situation.
\begin{lemma}\label{LemM3} 
If
 $\alpha,\beta{\in} L^2_\per(\Box)$, $\langle \alpha\beta\rangle{=}0$, $\alpha_\e(x){=}\alpha(x/\e)$, $\beta_\e(x){=}\beta(x/\e)$, 
  $\varphi,\psi\in H^1(\rd)$, then
  \begin{equation}\label{m.9}
( \alpha_\e S^\e\varphi,\beta_\e S^\e\psi) \le C\e^2
\langle \alpha^2\rangle^{1/2}\langle \beta^2\rangle^{1/2}
\|\nab\varphi\|\,\|\nab\psi\|,\quad C=const(d).
\end{equation}
\end{lemma}

 We now weaken 
 the conditions on the periodic function in  Lemma \ref{LemM2}. 
\begin{lemma}\label{LemM4} If
$\alpha,\beta{\in} L^2_\per(\Box)$, 
 $\alpha_\e(x){=}\alpha(x/\e)$, $\beta_\e(x){=}\beta(x/\e)$, 
  $\varphi{\in} L^2(\rd)$, and $\psi{\in} H^1(\rd)$, then
 \begin{equation}\label{m.10}
|( \alpha_\e S^\e\varphi,\beta_\e S^\e \psi) - \langle \alpha\beta\rangle
( \varphi,\psi)
|\le C\e
\langle \alpha^2\rangle^{1/2}\langle \beta^2\rangle^{1/2}
\|\varphi\|\,\|\nab\psi\|,\quad C=const(d).
\end{equation}
\end{lemma}
 The proof of (\ref{m.9}) and (\ref{m.10}) can be found in 
 \cite{P20a}--\cite{P20s}.

\textbf{4.2. Smoothing 
 with an  arbitrary kernel.} We consider the  smoothing operator
 \begin{equation}\label{m.11}
\Theta^\e\varphi(x)=\ild\varphi(x-\e\omega)\theta(\omega)\,d\omega, 
\end{equation}
where $\theta\in L^\infty(\rd)$ 
is compactly supported, $\theta\ge 0$, and $\int_{\rd} \theta(x)dx{=}1$. 

The estimates  (\ref{m.2})--(\ref{m.4}) formulated for the Steklov   smoothing operator remain valid for the general smoothig operator
 (\ref{m.7}) with only one note that right-hand side contains constants depending not only on the dimension $d$, but also on the kernel $\theta$. If, in addition,  the kernel $\theta$ is even, 
then $\Theta^\e$ possesses the properties of  type (\ref{m.6}) and (\ref{m.7}).
The following properties of the operator (\ref{m.11}) are highlighted in    \cite{P20s} and \cite{NY}. 

\begin{lemma}\label{LemM5} Assume that $\theta$ is piecewise $C^k$-smooth, $k$ is a natural number, 
$b\in L^2_\per(Y)$, $b_\e(x)=b(x/\e)$, and  
  $\varphi\in L^2(\rd)$. Then
   \begin{equation}\label{m.12}
\|\Theta^\e\nab^k\varphi\|\le C\e^{-k}
\|\varphi\|,\quad C=const(\theta),
\end{equation}
   \begin{equation}\label{m.13}
\|b_\e\Theta^\e\nab^k\varphi\|\le C\e^{-k}
\langle b^2\rangle^{1/2}\|\varphi\|,\quad C=const(\theta).
\end{equation}
\end{lemma}

\textbf{4.3. Steklov smoothing iteration.} It is obvious that the Steklov smoothing operator $S^\e$ can be defined by (\ref{m.11}) if the smoothing kernel  $\theta_1(x)$ is the characteristic function of the cube   $Y{=}[-{1}/{2},{1}/{2})^d$.
The double Steklov smoothing operator $(S^\e)^2=S^\e S^\e$ has the form (\ref{m.11}) with the smoothing kernel equal to the convolution  $\theta_2=\theta_1*\theta_1$.
Similarly, the triple Steklov smoothing operator $(S^\e)^3=S^\e {S}^\e S^\e$ is of the form (\ref{m.7}) with the smoothing kernel  $\theta_3=\theta_2*\theta_1$.
The kernels $\theta_2$ and $\theta_3$ (see computations for them in \cite{P20s}) turn to be piecewise $C^1$-smooth and
piecewise $C^2$-smooth respectively; as a consequence, the properties  (\ref{m.12}) and (\ref{m.13})  with $k=1$ hold surely
for the smoothing operator $\Theta^\e=(S^\e)^2$ or $\Theta^\e=(S^\e)^3$. In addition, since $\theta_2$ and $\theta_3$  are even,  the operators   $(S^\e)^2$ and $(S^\e)^3$ possess the properties of type (\ref{m.6}) and (\ref{m.7}).

Here is some corollary of Lemma \ref{LemM4}
involving iterations of  $S^\e$.

\begin{lemma}\label{LemM6} 
Under assumptions of Lemma \ref{LemM4} consider any generalized derivative $D_i\varphi$, $1\le i\le d$, of the function $\varphi{\in} L^2(\rd)$. Then 
 \begin{equation}\label{m.14}
|( \alpha_\e (S^\e)^3 D_i\varphi,\beta_\e S^\e \psi) - \langle \alpha\beta\rangle
(  (S^\e)^2D_i\varphi,\psi)
|\le C
\langle \alpha^2\rangle^{1/2}\langle \beta^2\rangle^{1/2}
\|\varphi\|\,\|\nab\psi\|  
\end{equation}
with the constant $C=const(d)$.
\end{lemma}
\doc
By Lemma \ref{LemM5}, both $L^2$-forms in left-hand side of (\ref{m.14}) make sense and are of order  $O(\e^{-1})$ as $\e\to 0$ because the smoothing kernels of $(S^\e)^3$
and $(S^\e)^2$ are, at least, piecewise $C^1$-smooth; therefore,  in both cases 
the generalized derivative $D_i\varphi$, smoothed by $(S^\e)^3$
or $(S^\e)^2$,
belongs to $\ld$  with the estimate of its 
norm due to (\ref{m.12}). Applying Lemma \ref{LemM4} to the pair of functions $\Phi=\e  (S^\e)^2D_i\varphi$ and $\psi$, we can write
\[
|( \alpha_\e S^\e\Phi,\beta_\e S^\e \psi) - \langle \alpha\beta\rangle
(  \Phi,\psi)
|\le C\e
\langle \alpha^2\rangle^{1/2}\langle \beta^2\rangle^{1/2}
\|\Phi\|\,\|\nab\psi\|,
\]
where $\|\Phi\|{=}\|\e  (S^\e)^2D_i\varphi\|{\le} C\|\varphi\|$ by Lemma \ref{LemM5}. Substituting here the expression for $\Phi$, we arrive at
\[
|\e( \alpha_\e (S^\e)^3D_i\varphi,\beta_\e S^\e \psi) - \e\langle \alpha\beta\rangle
(  (S^\e)^2D_i\varphi,\psi)
|\le C\e
\langle \alpha^2\rangle^{1/2}\langle \beta^2\rangle^{1/2}
\|\varphi\|\,\|\nab\psi\|,
\]
which after division by $\e$ gives
 (\ref{m.14}). The proof is completed.

We now specify 
the estimate (\ref{m.6}) for functions $\varphi\in H^4(\R^d)$
 and extend it at the same time to the case of an arbitrary  iteration $(S^\e)^k$.
 
\begin{lemma}\label{LemM7} 
 For any iteration of Steklov`s smoothing operator $(S^\e)^k$, $k{=}1,2,\ldots,$ the following inequality holds:
 \begin{equation}\label{m.15}
\|(S^\e)^k\varphi-\varphi-\e^2\gamma_{k}\,\Delta\,\varphi\|\le C\e^4\|\nab^4\varphi\|\quad \forall\varphi\in H^4(\R^d),\quad C=const(d,k), 
\end{equation}
where $\Delta$ is the Laplacean operator, i.e. $\Delta\,\varphi{=}D_jD_j\varphi$, the coefficient  $\gamma_k$ is determined by the smoothing kernel  $\theta_{k}$ of  $(S^\e)^k$ (see (\ref{m.17})), 
in particular,  the simple calcilation gives
$\gamma_{1}=1/24$, $\gamma_2=1/12$, $\gamma_3=1/8$.
\end{lemma}
\doc
We write the Taylor expansion with fourth-order remainder term
\[
\varphi(x+h)-\varphi(x)-D_i\varphi(x)h_i-\frac12 D_iD_j\varphi(x)h_ih_j
-\frac16 D_iD_jD_l\varphi(x)h_ih_jh_l=r_4(x,h),
\]
\[
r_4(x,h)=\frac16 \int_0^1 (1-t)^3 D_iD_jD_lD_m\varphi(x+th)h_ih_jh_lh_m\,dt.
\]
Substituting here $h{=}{-}\e\omega$ with $\omega{\in} Y{=}[-1/2,1/2)^d$, and thereafter multiplying both parts of the above equality by the smoothing kernel $\theta_{k}(\omega)$ of the operator $(S^\e)^k$, we integrate over $\omega\in Y$. As a result,  taking into account that
the function $\theta_{k}(\omega)$ is even, we obtain 
 \begin{equation}\label{m.16}
(S^\e)^k\varphi(x)-\varphi(x)-\e^2\gamma_{k}\,\Delta\,\varphi(x)=\e^4\int_Y r_4(x,\omega)\theta_{k}(\omega)\,d\omega
\end{equation}
with the coefficient
 \begin{equation}\label{m.17}
\gamma_{k}=\frac12 \int_Y \omega_1^2\theta_{k}(\omega)\,d\omega .
\end{equation}
To find the value of the integral (\ref{m.17}),
one can utilize 
the property of  separation of variables: $\theta_k(\omega)=
\chi_k(\omega_1)\cdots\chi_k(\omega_d)$; moreover,
functions
$\chi_k(t)$, 
 $t{\in}\mathbb{R}$, are restored successively one after another,
starting with $\chi_1(t)$, for all natural $k$
 (see more about 
 the smoothing kernels for iterations of $S^\e$ in Remark  4.8). By the Cauchy–-Schwarz inequality,  (\ref{m.16}) yields
\[
|(S^\e)^k\varphi(x)-\varphi(x)-\e^2\gamma_{k}\,\Delta\,\varphi(x)|^2\le \e^8 C\, \int_Y \int_0^1
|\nab^4\varphi(x-\e t\omega)|^2\,d\omega\,dt,
\]
and after integration with respect to $x$ estimate  
(\ref{m.15}) follows. 

The description of the smoothing kernels $\theta_1$, $\theta_2$, and $\theta_3$ given in Remark 4.8 allows us to compute the values of
$\gamma_{1}$, $\gamma_2$, and $\gamma_3$.

\bigskip
\noindent \textbf{Remark 4.8.} The approximation 
(\ref{2.13})--(\ref{2.15}) and Lemma \ref{LemM6} involve the iterations $(S^\e)^2$  and $(S^\e)^3$,
that is why we recall the expressions for their smoothing kernels, which were computed in \cite{P20s}.
The iteration $(S^\e)^2$ is the operator (\ref{m.11}) with the smoothing kernel
$\theta_2(x)=\chi_2(x_1)\ldots \chi_2(x_d)$, where $\chi_2(s)=\chi_1*\chi_1(s)$ and $\chi_1(s)$ is the characteristic function of the interval $[-\f{1}{2},\f{1}{2})$. Calculating the  convolution yields 
\[
\chi_2(s)=\begin{cases}
1-|s|,&\text{if $|s|\le 1$,}\\ 
0,&\text{if $|s|> 1$.}
\end{cases}
\]
Similarly,  $(S^\e)^3$ has the smoothing kernel
$\theta_3(x)=\chi_3(x_1)\ldots \chi_3(x_d)$, where $\chi_3(s)=\chi_1*\chi_2(s)$, $s\in \mathbb{R}$, and calculating the  convolution yields
\[
\chi_3(s)=\begin{cases}
3/4-s^2,&\text{if $|s|\le 1/2$,}\\ 
1/2(|s|-3/2)^2,&\text{if $ 1/2<|s|< 3/2$,}\\
0,&\text{if $|s|\ge 3/2$.}
\end{cases}
\]

\section{$\!\!\!\!\!\!$.  $L^2$-estimate with correctors }\label{sec5} 
 \setcounter{theorem}{0} \setcounter{equation}{0}

In this Section, we prove 
Theorem \ref{th5.1}. 

\textbf{5.1. Preliminaries.} From now on we denote by
 $b_\e$,
or $(b)_\e$, 
 the function that is obtained
 from
 periodic function $b(y)$ with periodicity cell $Y{=}[-\frac12 ,\frac12)^d$, 
 if we substitute $y{=}x/\e$, that is,
 \begin{equation} \label{5.1}
b_\e(x)=b(x/\e).
\end{equation}
For instance, $a_\e(x){=}a(x/\e)$, $N^{ij}_\e(x){=}N^{ij}(x/\e)$, $\left(a D N^{ij}\right)_\e{=}
\left(a(y) D_y N^{ij}(y)\right)|_{y=(x/\e)}$ and so on.

For a two-scale function, the differentiation formula is valid
 \begin{equation}\label{5.2}
 \ds{
 D\Phi(x,x/\varepsilon)
 } \atop\ds{
 =[D_x\Phi(x,y){+}\varepsilon^{-2}D_y\Phi(x,y){+}\varepsilon^{-1}(\n_x{\times}\n_y)\Phi(x,y)
 {+}\varepsilon^{-1}(\n_y{\times}\n_x)\Phi(x,y)
 ]|_{y=x/\varepsilon},
  }
\end{equation}
 where the differential operator 
   $(\n_x{\times}\n_y)$ sends a scalar function 
$\Phi(x,y)$ to the matrix of mixed second order derivatives  
$\{ \f{\partial^2  \Phi}{\partial x_i\partial y_j} \}_{i,j}$.
Using (\ref{5.2}), we have 
\begin{equation}\label{5.3}
 \ds{D \left(u(x){+}\e^2 N^{ij}(x/\e)D_{ij} u(x)\right)
 =Du(x){+}D_y N^{ij}(y)D_{ij} u(x){+} \varepsilon(\n_yN^{ij}(y))\times(\n D_{ij} u(x))
 }
 \atop\ds{+\varepsilon (\n D_{ij} u(x)) \times(\n_yN^{ij}(y))
+\varepsilon^2N^{ij}(y)D\,D_{ij} u(x), \quad 
 y=x/\varepsilon,
 }
\end{equation}
and from the symmetry property of the tensor $a(y)$ it follows that
\begin{equation}\label{5.4}
\ds{
a(x/\e)D \left(u(x)+\e^2 N^{ij}(x/\e)D_{ij} u(x)\right)=a(y)(D_yN^{ij}(y)+e^{ij})D_{ij}  u(x)}
\atop\ds{
{+}
2\varepsilon  a(y)(\n_yN^{ij}(y)){\times}(\n D_{ij}  u(x)){+}\varepsilon^2 a(y) N^{ij}(y)DD_{ij}  u(x),\,  y{=}x/\varepsilon,}
\end{equation}
where the matrix  $e^{ij}$ was introduced in the setting of problem (\ref{3.1}) on the periodicity cell. 
The calculations of the type (\ref{5.3}) and (\ref{5.4}) will be done several times in the following. 

\textbf{5.2. $H^2$-approximation.} The  pivotal role in the proof of
 estimate
 (\ref{2.17})  
 is played by 
 \begin{lemma}\label{Lem5.1} 
(i) Assume that the function $\tilde{u}^\e=u^{,\e}+\e^2U^\e_2+\e^3U^\e_3$ is specified  in (\ref{2.13})--(\ref{2.16}). Then its  discrepancy with  equation (\ref{2.1}) can be represented in the form
 \begin{equation} \label{5.5}
(A_\e+1)\tilde{u}^\e-f= D^*r_\e+r^0_\e+(f^{,\e}-f)=:F^\e,
\end{equation}
where
 \begin{equation} \label{5.6}
f^{,\e}=\Theta^\e f,\quad r^0_\e=\e^2U^\e_2+\e^3U^\e_3\stackrel{(\ref{2.14})}=
\e^2 N^{ij}_\e z_{ij}+\e^3 N^{ijk}_\e z_{ijk},
\end{equation}
 \begin{equation} \label{5.7}
  \ds{
 r_\e =
 \e^2 \left[2(\pa_q G^{ijk,pq})_\e z_{ijkp}+(aN^{ij})_\e Dz_{ij}+2 a_\e (\n N^{ijk}\times e^p)_\e 
z_{ijkp}+ G^{ij,kp}_\e z_{ijkp}\right]
}\atop\ds{
+\e^3  \left[ G^{ijk,pq}_\e z_{ijkpq}+(aN^{ijk})_\e Dz_{ijk}\right],
}
\end{equation}
 $\Theta^\e=(S^\e)^3$ is the triple Steklov smoothing operator,  and we use the notation
 \begin{equation} \label{5.8}
z_{ij}=D_iD_j u^{,\e},\quad z_{ijk}=D_iD_jD_k u^{,\e},\ldots, z_{ijkpq}=D_iD_jD_kD_pD_q u^{,\e}
\end{equation}
 for derivatives of the function $u^{,\e}=\Theta^\e \hat{u}^\e$ of order from second to fifth. All the above  functions
$ \hat{u}^\e$, $ N^{ij}$, $ N^{ijk}$, $ G^{ij,km}$, and $G^{ijk,pq}$ are defined in (\ref{2.10})
(\ref{3.1}), (\ref{3.9}), (\ref{3.8})
  and (\ref{3.13}) respectively.
  
(ii) For the right-hand side function $F^\e$ in (\ref{5.5}), we have the estimate
 \begin{equation} \label{5.9}
\|F^\e\|_{H^{-2}(\ld)}\le C\e^2\|f\|_{L^2(\rd)},\quad C=const (\lambda,d);
\end{equation}
 furthermore,    the term $r_\e$ 
 specified in (\ref{5.7}) satisfies the estimate
 \begin{equation} \label{5.90}
\|r^\e\|_{L^2(\rd)}\le C\e^2\|f\|_{L^2(\rd)},\quad C=const (\lambda,d).
\end{equation}
\end{lemma}

In fact, this assertion was 
proved  
in \cite{P21}, 
where the similar, as in (\ref{2.13})--(\ref{2.16}), approximation was studied with the  only difference in the smoothing operator ($(S^\e)^2$ in place of $(S^\e)^3$), which does not cause serious changes in arguments of
  \cite{P21}. In our approach, the structure of the right-hand side in (\ref{5.5}) 
determines the structure of
$L^2$-approximation from (\ref{2.20}); therefore, for the sake of completeness, we give here the proof of Lemma  \ref{Lem5.1} carrying it out in Section 6. Therein, 
 we show also
how 
Lemma  \ref{Lem5.1} provides $\e^2$ order
$H^2$-estimate  (\ref{2.17}) for the error of the approximation (\ref{2.13})--(\ref{2.16}).
The latter  estimate itself is essentially used
in the proof of  
(\ref{2.20}). 

\textbf{5.3. $L^2$-approximation.} 
In what follows we denote  by $c$ or $C$ a general positive constant, depending only on the dimension $d$
and the ellipticity constant
$\lambda$
from (\ref{2.2}), possibly varying  from line to line.
We continue to use the simplified notation introduced in  Secton 4 for the norm and the inner product of $\ld$.

The proof of estimate 
(\ref{2.20}) is splitted into several steps.

1$^\circ$ From (\ref{2.17}), by mere weakening the norm, it follows  
 \begin{equation}\label{5.11}
 \|u^\e-\tilde{u}^\e\|\le c\e^2\|f\|. 
\end{equation}
Due to smoothing properties, the approximation in (\ref{5.11}) can be simplified by discarding both  correctors without upsetting the order of the estimate (see Remark 
2 in Section 7). Bur right now our goal is to improve the accuracy of the approximation for
$u^\e$ in $L^2$-norm uo to $\e^3$ order   by means of additional correctors. To find this kind of correctors, we evaluate more sharply 
the norm $\|u^\e-\tilde{u}^\e\|$,
studying the  $L^2$-form
\begin{equation}\label{5.12}
I:=(u^\e-\tilde{u}^\e,h), \quad h\in L^2(\rd).
\end{equation}
To this end, 
we introduce another version of equation (\ref{2.1}), with $h\in \ld$ on the right-hand side,
\begin{equation}\label{5.13}
v^\e\in H^2(\rd),\quad (A_\e+1)v^\e=h;
\end{equation}
 the corresponding homogenized equation
\begin{equation}\label{5.14}
(\hat{A}_\e+1)\hat{v}^\e=h;
\end{equation}
and  the $H^2$-approximation of the form
\beq\label{5.15}
\tilde{v}^\e=v^{,\e}(x)+\e^2\,V_2^\e(x)+\e^3\,V_3^\e(x),
\eeq
where
\beq\label{5.16}
V_2^\e(x)=N^{ij}(x/\e)D_iD_{j}v^{,\e}(x),\quad V_3^\e(x)=N^{ijk}(x/\e)D_iD_{j}D_kv^{,\e}(x),
\eeq
\beq\label{5.17}
v^{,\e}(x)=\Theta^\e\, \hat{v}^\e(x),\quad \Theta^\e=S^\e S^\e S^\e
\eeq
 with the functions
 $N^{ij}(y)$, $N^{ijk}(y)$, and $ \hat{v}^\e(x)$  defined before. The function $\tilde{v}^\e$ is quite similar  to the approximation (\ref{2.13})--(\ref{2.16}) of $\tilde{u}^\e$.
 For the difference of the solution to (\ref{5.13}) and the function $\tilde{v}^\e$, the 
estimate  of type (\ref{2.17}) holds, namely, 
\begin{equation}\label{5.18}
\|v^\e-\tilde{v}^\e\|_{H^2(\rd)}\le c\e^2\|h\|. 
\end{equation}
We recall also that the uniform in $\e$ estimate for the solution of the homogenized equation (\ref{5.14}) holds:
\begin{equation}\label{5.19}
\|\hat{v}^\e\|_{ H^{4}(\rd)}\le c\|h\|, 
\end{equation}
which is similar to (\ref{2.12}).

We start with the following representation
\[
I=(u^\e-\tilde{u}^\e,h)\stackrel{(\ref{5.13})}=(u^\e-\tilde{u}^\e, ({A}_\e+1){v}^\e)=
(({A}_\e+1)(u^\e-\tilde{u}^\e), {v}^\e)
\]
\[
\stackrel{(\ref{2.1})}=(f{-}({A}_\e+1)\tilde{u}^\e, {v}^\e)
\stackrel{(\ref{5.5})}=-(F^\e, {v}^\e)=-(F^\e, v^\e-\tilde{v}^\e)-(F^\e, \tilde{v}^\e),
\]
where
$
|(F^\e, v^\e-\tilde{v}^\e)|\le c \e^4\|f\|\,\|h\|
$
 in view of (\ref{5.9}) and (\ref{5.18}),  thereby yielding
\begin{equation}\label{5.21}
I\cong -(F^\e, \tilde{v}^\e).
\end{equation}
Here and below, we let $\cong$ denote 
the approximate equality obtained from the exact equality by discarding terms $T$ admitting an
estimate
\[
|T|\le c\e^3\|f\|\,\|h\|,\quad c=const(d,\lambda).
\]
An expression 
 $T$ satisfying this estimate is said to be \textit{nonessential}.

In view of (\ref{5.5}), we can rewrite (\ref{5.21}) as
\begin{equation}\label{5.22}
I\cong -(D^* r^\e, \tilde{v}^\e)
-(r_0^\e, \tilde{v}^\e)
+(f-f^{,\e}, \tilde{v}^\e):=I_1+I_2+I_3,
\end{equation}
where each form 
 $I_i$ should be estimated.

2$^\circ$ 
 Starting with the representation
\[
I_2:=-(r_0^\e,  \tilde{v}^\e)\stackrel{(\ref{5.6})}=
-(\e^2 N^{ij}_\e z_{ij}+\e^3 N^{ijk}_\e z_{ijk}, \tilde{v}^\e)=
-\e^2( N^{ij}_\e z_{ij}, \tilde{v}^\e)-\e^3 (N^{ijk}_\e z_{ijk}, \tilde{v}^\e)
,
\]
we are going  to show that the both terms on the right-hand side are nonessential.
By Lemma \ref{LemM1} (in its version with the smoothing operator $\Theta^\e=(S^\e)^3$), we get
$$\|N^{ijk}_\e z_{ijk}\|\stackrel{(\ref{5.8})}=\|N^{ijk}_\e D_iD_j D_k\Theta^\e \hat{u}^\e\|\stackrel{(\ref{m.4})}\le c\langle |N^{ijk}|^2\rangle^{1/2}\|\nab^3\hat{u}^\e\|\stackrel{(\ref{2.12})}
\le C \|f\|.
$$
Similarly, by Lemma \ref{LemM1} combined with (\ref{5.19}), we have
 \begin{equation} \label{5.23}
\|\tilde{v}^\e\|\stackrel{(\ref{5.14})}=\|v^{,\e}{+}\e^2\,V_2^\e{+}\e^3\,V_3^\e\|
\stackrel{(\ref{5.16})}\le 
\|v^{,\e}\|{+}\e^2\|N^{ij}_\e D_iD_{j}v^{,\e}\|{+}\e^3\|N^{ijk}_\e D_iD_{j}D_kv^{,\e}\|{\le} c\|h\|,
\end{equation}
where we use also that $v^{,\e}=\Theta^\e\hat{v}^\e$.
From the last two estimates, it is clear  that
\[
\e^3 (N^{ijk}_\e z_{ijk}, \tilde{v}^\e)
\le \e^3 \|N^{ijk}_\e z_{ijk}\|\,\| \tilde{v}^\e\|\cong 0.
\]
On the other hand, by Lemma  \ref{LemM2}, 
\[
\e^2 (N^{ij}_\e z_{ij}, \tilde{v}^\e)=\e^2 (N^{ij}_\e S^\e \varphi, \tilde{v}^\e)
\stackrel{(\ref{m.5})}\le C\e^3\langle |N^{ij}|^2\rangle^{1/2} \|\varphi\|\,\| \nab\tilde{v}^\e\|\cong 0.
\]
Here, $z_{ij}=S^\e \varphi$  with $\varphi=(S^\e)^2 D_iD_{j}\hat{u}^\e$ due to (\ref{5.8}); by property of  solutions of the problem (\ref{3.1}), 
 we have $\langle N^{ij}\rangle=0$; and
 we know that
$\|\varphi\|\le c\|f\| $  in view of (\ref{2.12}) and $\|\nab\tilde{v}^\e\|\le c\|h\|$ similarly, as (\ref{5.23}). Ultimately, 
\begin{equation}\label{5.24}
I_2\cong 0.
\end{equation}

3$^\circ$ 
The first form in (\ref{5.22}) can be simplified as
\[ 
I_1:=-(D^* r^\e, \tilde{v}^\e)=-( r^\e,D \tilde{v}^\e)\cong
-( r^\e,(D N^{mn}+e^{mn})_\e w_{mn}),
\]
because estimate (\ref{5.90}) is valid for the term $ r^\e$ and calculations of type (\ref{5.3}), (\ref{5.4}) give here
\[ 
D \tilde{v}^\e\stackrel{(\ref{5.15})-(\ref{5.17})}=(D N^{mn}+e^{mn})_\e w_{mn}+O(\e),
\]
where  the term  $O(\e)$ includes all the summands whose $\ld$-norms are less than $c\e\|h\|$.
We use the notation
 (similar to (\ref{5.8}))
 \begin{equation} \label{5.25}
w_{mn}=D_mD_n v^{,\e}, 
\ldots, w_{mnkpq}=D_mD_nD_kD_pD_q v^{,\e}
\end{equation}
for the derivatives of $v^{,\e}=\Theta^\e \hat{v}^\e$ of orders from second to fifth. 

Taking  into account (\ref{5.7}), we rewrite $I_1$ as a sum  
\begin{equation} \label{5.26}
\ds{ 
I_1\cong
-( \e^2 [2(\pa_q G^{ijk,pq})_\e z_{ijkp}+(aN^{ij})_\e Dz_{ij}+2 a_\e (\n N^{ijk}\times e^p)_\e 
z_{ijkp}+ G^{ij,kp}_\e z_{ijkp}]
}\atop\ds{
+\e^3   G^{ijk,pq}_\e z_{ijkpq}+\e^3(aN^{ijk})_\e Dz_{ijk}
,(D N^{mn}+e^{mn})_\e w_{mn}):=I_{11}+\ldots+I_{16}.}
\end{equation} 
 Each of  the forms 
 $I_{11},\ldots,I_{16}$ should be estimated.
 
By Lemma \ref{LemM4},
\[
I_{11}{:=}-2\e^2 ( (\pa_q G^{ijk,pq})_\e z_{ijkp},(D N^{mn}{+}e^{mn})_\e w_{mn}){\cong}
-2\e^2 \langle\pa_q G^{ijk,pq}{\cdot}( D N^{mn}{+}e^{mn})\rangle ( \varphi_{ijkp}, \psi_{mn}).
\]
To explain this equality, first, note that
$
z_{ijkp}{=}S^\e\varphi_{ijkp}$, $\varphi_{ijkp}{=}
D_iD_jD_kD_p (S^\e)^2 \hat{u}^\e$ with $\varphi_{ijkp}\in \ld$ such that
$
\|\varphi_{ijkp}\|\le c\|f\|
$, by properties of smoothing and estimate  (\ref{2.12}). In a similar way,
 according to (\ref{5.25}),
$w_{mn}{=}S^\e\psi_{mn}$, $\psi_{mn}{=}D_mD_n (S^\e)^2 \hat{v}^\e$ and $\psi_{mn}\in H^1(\rd)$
such that 
$
\|\nab \psi_{mn}\|\le c\|h\|,
$
by (\ref{5.19}). Therefore,  writing the latter representation of $I_{11}$, we   discard its part that is
nonessential, due to (\ref{m.10}). Since
 \[
 \langle\pa_q G^{ijk,pq}{\cdot}( D N^{mn}{+}e^{mn})\rangle=
 \langle\pa_q G^{ijk,pq}{\cdot} D N^{mn}\rangle=
 -\langle D^* G^{ijk,pq}  \pa_q N^{mn}\rangle
\stackrel{(\ref{3.13})_1}=- \langle g^{ijk}_{pq}  \pa_q N^{mn}\rangle,
 \]
we finally obtain
  \begin{equation} \label{5.27}
I_{11}\cong 2\e^2 \langle g^{ijk}_{pq}  \pa_q N^{mn}\rangle ( \varphi_{ijkp}, \psi_{mn}).
\end{equation} 

To analyse the next summand in (\ref{5.26}), we use first the symmetry property of the tensor $a$:
\[
I_{12}{:=}-\e^2 (  (aN^{ij})_\e Dz_{ij}, (D N^{mn}+e^{mn})_\e w_{mn})=-\e^2 (  N^{ij}_\e Dz_{ij}, 
a_\e(D N^{mn}+e^{mn})_\e w_{mn}),
\]
and then recall that $a(D N^{mn}+e^{mn})=g^{mn}+\hat{a}e^{mn}$  in view of (\ref{3.6}). 
By Lemma \ref{LemM2},
\[
\e^2 (  N^{ij}_\e Dz_{ij}, 
\hat{a}e^{mn} w_{mn})\cong 0,
\]
because $\langle N^{ij}\rangle=0$, $Dz_{ij}=S^\e\varphi$ with $\|\varphi\|\le c
\|f\|$  due to (\ref{5.8}) and (\ref{2.12}), and also 
$\hat{a}e^{mn} w_{mn}\in H^1(\rd)$ such that
$\|\nab(\hat{a}e^{mn} w_{mn})\|\le c\|h\|$ in view of (\ref{5.25}) and (\ref{5.19}).
As a result, discarding nonessential terms, we obtain 
 \begin{equation} \label{5.28}
I_{12}{\cong}{-}\e^2 (  N^{ij}_\e Dz_{ij}, g^{mn}_\e w_{mn})
{=}{-}\e^2 (   z_{ijkp}, (N^{ij}g^{mn}_{kp})_\e w_{mn})\stackrel{(\ref{m.10})}
\cong-\e^2 \langle N^{ij}g^{mn}_{kp}\rangle(   \varphi_{ijkp},  \psi_{mn}),
\end{equation}
where we write  $Dz_{ij}\cdot g^{mn}_\e=
z_{ijkp}(g^{mn}_{kp})_\e$ according to the definition of the operator $D$ combined  with the convention  (\ref{5.8}). Besides, we introduce in this write-up
the new notation (similar to (\ref{5.8}) and (\ref{5.25}))
 \begin{equation} \label{5.29}
\varphi_{mn}=D_mD_n (S^\e)^2\hat{u}^\e, 
\ldots, \varphi_{mnkpq}=D_mD_nD_kD_pD_q (S^\e)^2\hat{u}^\e
\end{equation}
for the derivatives of the smoothed function $(S^\e)^2\hat{u}^\e$  of order from second to fifth, and likewise
 \begin{equation} \label{5.30}
\psi_{mn}=D_mD_n (S^\e)^2\hat{v}^\e, 
\ldots, \psi_{mnkpq}=D_mD_nD_kD_pD_q (S^\e)^2\hat{v}^\e
\end{equation}
for the derivatives of the smoothed function $(S^\e)^2\hat{v}^\e$  of order from second to fifth.

We simplify the next two forms   $I_{1i}$
in (\ref{5.26}) in a similar fashion:
\[
I_{13}{:=}-2\e^2 (  a_\e(\nab N^{ijk}\times e^p)_\e z_{ijkp}, (D N^{mn}+e^{mn})_\e w_{mn})\]
\[
{=}-2\e^2 ( (\nab N^{ijk}\times e^p)_\e z_{ijkp}, 
a_\e(D N^{mn}+e^{mn})_\e w_{mn})
{=}-2\e^2 ( (\nab N^{ijk}\times e^p)_\e z_{ijkp}, 
(g^{mn}_\e+\hat{a}e^{mn}) w_{mn})
\]
\[ 
\stackrel{(\ref{m.5})} \cong -2\e^2 ( (\nab N^{ijk}\times e^p)_\e z_{ijkp}, 
g^{mn}_\e w_{mn})\stackrel{(\ref{m.10})}\cong -2\e^2 \langle (\nab N^{ijk}\times e^p)\cdot g^{mn}\rangle(  \varphi_{ijkp}, 
 \psi_{mn})
,
\]
or more shortly,
\begin{equation} \label{5.31}
I_{13}\cong -2\e^2 \langle\partial_q N^{ijk}g^{mn}_{qp}\rangle(  \varphi_{ijkp}, 
 \psi_{mn}),
\end{equation}
where we take into account the symmetry of the tensor $a$, the definition of the matrix $g^{mn}$, the  constancy of the matrix $\hat{a}e^{mn}$, and
the equality $\langle \n N^{ijk}\rangle=0$, combined with conventions (\ref{5.29}), (\ref{5.30}) and estimates (\ref{2.12}),
(\ref{5.19}). Likewise,
 \begin{equation} \label{5.32}
 \ds{
I_{14}{:=}-\e^2 (   G^{ij,kp}_\e z_{ijkp}, (D N^{mn}+e^{mn})_\e w_{mn})\stackrel{(\ref{m.5})}\cong
-\e^2 (   G^{ij,kp}_\e z_{ijkp}, (D N^{mn})_\e w_{mn})
}\atop\ds{
\stackrel{(\ref{m.10})}\cong-\e^2
\langle G^{ij,kp}\cdot D N^{mn}\rangle(  \varphi_{ijkp},  \psi_{mn})
\stackrel{(\ref{3.13})_1}=-\e^2
\langle g^{ij}_{kp}N^{mn}\rangle(  \varphi_{ijkp}, 
 \psi_{mn}),}
\end{equation}
where,  at the first step we show
 \[
\e^2 (   G^{ij,kp}_\e z_{ijkp}, e^{mn} w_{mn})\cong 0,
\]
by  Lemma \ref{LemM2}, thanks to the property $\langle G^{ij,kp}\rangle=0$.

To estimate the last two forms in (\ref{5.26}), we use  some new arguments. Firstly, we discard terms
 that are nonessential  by Lemma \ref{LemM6}:
\begin{equation} \label{5.33}
I_{15}{:=}-\e^3 (   G^{ijk,pq}_\e z_{ijkpq}, (D N^{mn}+e^{mn})_\e w_{mn})
\stackrel{(\ref{m.14})}\cong-\e^3
\langle  g^{ijk}_{pq} N^{mn}\rangle(  \varphi_{ijkpq},  \psi_{mn}),
\end{equation}
where we have passed to the notation
  (\ref{5.29}), (\ref{5.30}) and used the equality
\[
\langle  G^{ijk,pq}\cdot( D N^{mn}+e^{mn})\rangle=
\langle  G^{ijk,pq}\cdot D N^{mn}\rangle\stackrel{(\ref{3.13})_1}=
\langle  g^{ijk}_{pq} N^{mn}\rangle.
\]
But the 
right-hand side of (\ref{5.33}) itself can be regarded as nonessential, because integration by parts ensures the required bound: 
\begin{equation} \label{5.34}\ds{
(  \varphi_{ijkpq},  \psi_{mn})\stackrel{(\ref{5.29}), (\ref{5.30})}=
(D_iD_jD_kD_pD_q (S^\e)^2\hat{u}^\e,
D_mD_n (S^\e)^2\hat{v}^\e)
}\atop\ds{
=-(D_jD_kD_pD_q (S^\e)^2\hat{u}^\e,
D_iD_mD_n (S^\e)^2\hat{v}^\e)\le \|\nab^4 \hat{u}^\e\|\,\|\nab^3 \hat{v}^\e\|
\stackrel{(\ref{2.12}), (\ref{5.19})}\le c\|f\|\,\|h\|.}
\end{equation}
Consequently, 
 \begin{equation} \label{5.35}
I_{15} \cong 0.  
\end{equation}
Likewise, 
\[
I_{16}{:=}-\e^3 (  (aN^{ijk})_\e Dz_{ijk}, (D N^{mn}+e^{mn})_\e w_{mn})
=-\e^3(  N^{ijk}_\e Dz_{ijk}, a_\e(D N^{mn}+e^{mn})_\e w_{mn})
\]
\[
\stackrel{(\ref{3.6})}=-\e^3(  N^{ijk}_\e Dz_{ijk},(g_\e^{mn}+\hat{a} e^{mn}) w_{mn})
\stackrel{(\ref{m.14})}\cong-\e^3 (  D\varphi_{ijk},
\langle  N^{ijk}g^{mn}\rangle 
 \psi_{mn});
\]
therefore,  because of the constancy of the matrix $\langle  N^{ijk}g^{mn}\rangle$,  estimate (\ref{5.34}) yields
 \begin{equation} \label{5.36}
I_{16} \cong  0. 
\end{equation}

From (\ref{5.26})--(\ref{5.36}) we obtain
  \begin{equation} \label{5.37}\ds{
  I_1\cong I_{11}+\ldots+I_{14}}
  \atop\ds{\cong
\e^2 \left(2\langle g^{ijk}_{pq}  \pa_q N^{mn}\rangle-2\langle\partial_q N^{ijk}g^{mn}_{qp}\rangle-\langle N^{ij}g^{mn}_{kp}\rangle-
\langle g^{ij}_{kp}N^{mn}\rangle\right) ( \varphi_{ijkp}, \psi_{mn}).}
\end{equation}

4$^\circ$ Now, we estimate the third form in (\ref{5.22}): 
 \begin{equation} \label{5.38}
\ds{
I_3:=(f-f^{,\e}, \tilde{v}^\e)\stackrel{(\ref{5.15})}=
(f-f^{,\e},v^{,\e}+\e^2\,V_2^\e+\e^3\,V_3^\e)}
\atop\ds{
{=}
(f{-}f^{,\e},v^{,\e}){+}(f,\e^2\,V_2^\e{+}\e^3\,V_3^\e){-}(f^{,\e},\e^2\,V_2^\e{+}\e^3\,V_3^\e)\cong
(f{-}f^{,\e},v^{,\e}){+}(f,\e^2\,V_2^\e),}
\end{equation}
where we    omit some  terms that are nonessential by lemmas \ref{LemM1},  which is applied  to the two forms containing the factor  $
e^3$,  and by Lemma \ref{LemM2}, which is applied to the form
$$
(f^{,\e},\e^2\, V_2^\e)\stackrel{(\ref{5.6})}= \e^2 (f^{,\e}, N^{ij}_\e z_{ij})\cong 0.
$$
   In addition, recalling that $\Theta^\e{=}(S^\e)^3$, 
 we study the term
   \[
  (f-f^{,\e},v^{,\e})\stackrel{(\ref{5.6})}=(f-\Theta^\e f,v^{,\e})=(f,v^{,\e}-\Theta^\e v^{,\e}){=}
  (f,v^{,\e}{+}\e^2\gamma_{3}\,\Delta\,v^{,\e}{-}\Theta^\e v^{,\e}){-}(f,\e^2\gamma_{3}\,\Delta\,v^{,\e}),
  \]
  using Lemma \ref{LemM7}. Here,
    \[
    (f,v^{,\e}{+}\e^2\gamma_{3}\,\Delta\,v^{,\e}{-}\Theta^\e v^{,\e})\le \|f\|\,\|v^{,\e}{+}\e^2\gamma_{3}\,\Delta\,v^{,\e}{-}\Theta^\e v^{,\e}\|\stackrel{(\ref{m.15})}\le C\e^4 \|f\|\|v^{,\e}\|
  \]
and
\[
\|v^{,\e}\|\stackrel{(\ref{5.17})}\le c\|\hat{v}^\e\|
\stackrel{(\ref{5.19})}\le C\|h\|,
\]
whence
$(f-f^{,\e},v^{,\e})\cong -(f,\e^2\gamma_{3}\,\Delta\,v^{,\e}) $. 
Finally,  (\ref{5.38}) is rewritten as
 \begin{equation} \label{5.39}
I_3\cong(f,\e^2\,V_2^\e-\e^2\gamma_{3}\,\Delta\,v^{,\e})\cong
(f,\e^2\,V_2^\e-\e^2\gamma_{3}\,\Delta\,\hat{v}^{\e}),
\end{equation}
where at the last step the function $v^{,\e}=\Theta^\e \hat{v}^{\e}$ is replaced with
$\hat{v}^{\e}$  thanks to the estimate
\[
\|\Theta^\e\Delta \hat{v}^{\e}-\Delta\hat{v}^{\e}\|\stackrel{(\ref{m.6})}\le C\e^2 \|\n^4\hat{v}^{\e}\|\stackrel{(\ref{2.12})}\le C\e^2 \|h\|.
\]

5$^\circ$ From (\ref{5.12}) anf (\ref{5.22}) combined with (\ref{5.24}), (\ref{5.37}), and (\ref{5.39}), we get
\begin{equation}\label{5.40}
(u^\e-\tilde{u}^\e,h)\cong (f,\e^2\,V_2^\e-\e^2\gamma_{3}\,\Delta\,\hat{v}^{\e} +
\e^2 c_{ijkpmn}( \varphi_{ijkp}, \psi_{mn}), 
\end{equation}
where
\begin{equation}\label{5.41}%
c_{ijkpmn}=2\langle g^{ijk}_{pq}  \pa_q N^{mn}\rangle-2\langle\partial_q N^{ijk}g^{mn}_{qp}\rangle-\langle N^{ij}g^{mn}_{kp}\rangle-
\langle g^{ij}_{kp}N^{mn}\rangle.
\end{equation}
We make now some 
transformation in (\ref{5.40}). Instead of the function
$\tilde{u}^\e=u^{,\e}+\e^2\,U_2^\e+\e^3\,U_3^\e$ defined in (\ref{2.13})--(\ref{2.16}), we leave only  the sum of 
the first two terms $u^{,\e}+\e^2\,U_2^\e$  discarding the term with the factor $\e^3$ which yields the form 
$(\e^3\,U_3^\e,h)\cong 0$ by properties of smoothing. We recall that
$U_2^\e=N_\e\cdot Du^{,\e}$ and
$V_2^\e=N_\e\cdot Dv^{,\e}$  in view of (\ref{2.14}) and (\ref{5.16})  (where for brevity the matrix $N=\{N^{ij}\}_{i,j}$, composed of the solutions $N^{ij}$ to the problem (\ref{3.1}),
is used); besides, according to (\ref{2.15})
 and (\ref{5.17}), $u^{,\e}=\Theta^\e \hat{u}^\e$ and
$v^{,\e}=\Theta^\e \hat{v}^\e$, where $\Theta^\e=(S^\e)^3$.
The forms in (\ref{5.40}), which contain the correctors
 $U_2^\e=N_\e\cdot(S^\e)^3 D\hat{u}^\e$ and
$V_2^\e=N_\e\cdot(S^\e)^3 D\hat{v}^\e$,
can be simplified in the smoothing operator. For instance,
\[
(\e^2\,U_2^\e,h)=(\e^2N_\e\cdot(S^\e)^3D\hat{u}^\e,h)
=(\e^2N_\e\cdot((S^\e)^3\pm S^\e)D\hat{u}^\e,h)\cong (\e^2N_\e\cdot S^\e D\hat{u}^\e,h),
\]
since Lemma \ref{LemM1} and smoothing properties yield the estimate
\[
\|N_\e\cdot((S^\e)^3- S^\e)D\hat{u}^\e\|\stackrel{(\ref{m.4})}\le \langle |N|^2\rangle^{1/2}
\|((S^\e)^2- 1)D\hat{u}^\e\|\stackrel{(\ref{m.6})}\le c\e^2\|\n^4\hat{u}^\e\|
\stackrel{(\ref{2.12})}\le C\e^2\|f\|.
\]
 Similarly, $(f,\e^2\,V_2^\e)\cong (f,\e^2 N_\e\cdot S^\e D\hat{v}^\e).
$
  As for smoothing in the form
$( \varphi_{ijkp}, \psi_{mn})$, it can be omitted at all. In fact, 
first,  expressions for $\varphi_{ijkp}$ and $\psi_{mn}$ from (\ref{5.29}) and (\ref{5.30}) are written; then, integrating by parts leads to
\[
\e^2 ( \varphi_{ijkp}, \psi_{mn})=-\e^2(D_jD_kD_p (S^\e)^2\hat{u}^\e,D_iD_mD_n (S^\e)^2\hat{v}^\e)\cong
-\e^2(D_jD_kD_p \hat{u}^\e,D_iD_mD_n \hat{v}^\e),
\]
where the property of the type (\ref{m.3}) for  $(S^\e)^2$ is used  at the last step.

We now summarize the above changes and substitute the expressions
$ \hat{u}^\e=(\hat{A}_\e+1)^{-1}f$,
$ \hat{v}^\e=(\hat{A}_\e+1)^{-1}h$, and $ u^\e=(A_\e+1)^{-1}f$  
into (\ref{5.40}), thus getting
\begin{equation}\label{5.42}
\left((A_\e+1)^{-1}f-\Theta^\e(\hat{A}_\e+1)^{-1}f-\e^2 (K_2(\e)+(K_2(\e))^*+M_\e+L_\e)f,h\right)\cong 0,
\end{equation}
where 
\begin{equation}\label{5.43}
\ds{
K_2(\e)=N_\e\cdot S^\e D\, (\hat{A}_\e+1)^{-1},\quad
M_\e=-\gamma_{3}\,\Delta\,(\hat{A}_\e+1)^{-1},
}\atop\ds{
L_\e=(\hat{A}_\e+1)^{-1}c_{ijkpmn}D_{ijkpmn} (\hat{A}_\e+1)^{-1},\quad
D_{ijkpmn}=D_iD_jD_kD_pD_mD_n,}
\end{equation}
$(K_2(\e))^*$ is the adjoint operator, $\Delta$ is the Laplacean operator, the coefficients $c_{ijkpmn}$ for all indices from 1 to d are determined in (\ref{5.41}), and $\gamma_{3}=1/8$ (see 
Lemma \ref{LemM7}). 
Due to
the property (\ref{m.15}) and the definition of 
$M_\e$ from (\ref{5.43}),
the second term in (\ref{5.42}) can be rewritten, in its turn,  as
 $$\Theta^\e(\hat{A}_\e+1)^{-1}f=(\hat{A}_\e+1)^{-1}f+\e^2\gamma_3\Delta (\hat{A}_\e+1)^{-1}f+O(\e^4)=
(\hat{A}_\e+1)^{-1}f-\e^2M_\e f+O(\e^4).
$$
Hence,  
according to the convention on $\cong$ we have adopted (see the explanation for (\ref{5.21})), it follows from (\ref{5.42}) that
\[
\|(A_\e+1)^{-1}f-(\hat{A}_\e+1)^{-1}f-\e^2 (K_2(\e)+(K_2(\e))^*+2M_\e+L_\e)f  \|
\le C\ve^3\|f\|,
\]
which means
(\ref{2.20}) with the function $w^\e=(\hat{A}_\e+1)^{-1}f+\e^2 (K_2(\e)+(K_2(\e))^*+2M_\e+L_\e)f.$ 

Eventually, we have established

\begin{theorem}\label{th5.1} Let $\hat{A}_\e$ be an operator from (\ref{2.8}), 
 and consider $\K(\e){=}K_2(\e){+}(K_2(\e))^*{+}2M_\e+L_\e$ with the operators $K_2(\e)$, $M_\e$,  $L_\e$ defined in (\ref{5.43}) and (\ref{5.41}). Then the sum 
 $(\hat{A}_\e{+}1)^{-1}+\e^2 \K(\e)$ approximates the resolvent $(A_\e{+}1)^{-1}$ in the operator 
 $L^2$-norm  
 in such a way that (\ref{2.22}) holds.
\end{theorem}

\section{$\!\!\!\!\!\!$. Approximation in $H^2$-norm}\label{sec6}
 \setcounter{theorem}{0} \setcounter{equation}{0}

\textbf{6.1.} 
We shall first provide a proof of Lemma \ref{Lem5.1}, and then show that estimate (\ref{2.17}) is a corollary of it.
The  discrepancy  of the function (\ref{2.13}) in  equation (\ref{2.1})
admits the representation
\[
A_\e \tilde{u}^\e{+}\tilde{u}^\e{-}f=(A_\e{+}1)\tilde{u}^\e{-}(\hat{A}_\e{+}1) u^{,\e}{+}(f^{,\e}{-}f)=
(A_\e \tilde{u}^\e{-}\hat{A}_\e u^{,\e}){+}
(v^\e{-} u^{,\e}){+}(f^{,\e}{-}f),
\]
where at the first step we take into account the equality
\[
(\hat{A}_\e+1) u^{,\e}=f^{,\e},\quad f^{,\e}=\Theta^\e f,
\]
obtained  by applying the smoothing operator $\Theta^\e$ to the both parts of (\ref{2.10})
 if the notation  $\Theta^\e\hat{u}^\e=u^{,\e}$ is used (see (\ref{2.15})).
Hence, 
 \begin{equation} \label{05.6}
 \ds{
A_\e \tilde{u}^\e{+}\tilde{u}^\e{-}f\stackrel{(\ref{2.13}), (\ref{2.1}), (\ref{2.9})}
=D^*\left(a_\e D(u^{,\e}+\e^2\,U_2^\e+\e^3\,U_3^\e)-\hat{a}Du^{,\e} -\e b^{rst}D_{rst}u^{,\e}\right)
}\atop\ds{
+(\e^2\,U_2^\e+\e^3\,U_3^\e)+(f^{,\e}-f).}
\end{equation}
To derive from here the representation (\ref{5.5}), we need to study the first term in the right-hand side of (\ref{05.6}), more exactly,  
the difference of flows 
 \begin{equation} \label{05.7}
R_\e:=a_\e D(u^{,\e}+\e^2\,U_2^\e+\e^3\,U_3^\e)-\hat{a}Du^{,\e} -\e b^{rst}D_{rst}u^{,\e}.
\end{equation}
Calculations similar to (\ref{5.3}) and (\ref{5.4}) combied with the notation (\ref{5.1}) and (\ref{5.8}) yield
\[
R_\e\stackrel{
(\ref{2.14})}
=a_\e[D u^{,\e}+z_{ij}(DN^{ij})_\e+2\e(\n N^{ij})_\e\times \n z_{ij}+\e^2N^{ij}_\e D z_{ij}
\]
\[
+\e z_{ijk}(DN^{ijk})_\e+2\e^2(\n N^{ijk})_\e\times \n z_{ijk}+\e^3 N^{ijk}_\e D z_{ijk}
]- \hat{a}Du^{,\e} -\e b^{ijk}z_{ijk}.
\]
Collecting the terms with the same factor
 $\e^n$, $n\ge 0$,
we get
\begin{equation} \label{05.9}
\ds{
R_\e=\left(a_\e(DN^{ij}+e^{ij})-\hat{a}e^{ij}\right)_\e z_{ij}+
\e\left[\left(a DN^{ijk}+2a(\n N^{ij}{\times} e^k)- b^{ijk}
\right)_\e z_{ijk}\right]
}\atop\ds{
 +\e^2\left[(aN^{ij})_\e Dz_{ij}+2 a_\e  (\n N^{ijk} {\times} e^m)_\e
z_{ijkm}\right]+
e^3\left[(aN^{ijk})_\e Dz_{ijk}\right]:=R_\e^0+R_\e^1+R_\e^2+R_\e^3,}
\end{equation}
where the vectors $e^m$ with   $m$ from 1 to d  form the canonical base in $\rd$.
  The term  $R_\e^0$ from the above sum contains the oscillating function
$
(a(DN^{ij}+e^{ij})-\hat{a}e^{ij})_\e
=g^{ij}_\e
$ (see (\ref{3.6})) and, thus, can be written
 shortly as $R_\e^0=g^{ij}_\e z_{ij}$, where the  1-periodic matrix $g^{ij}=\{g^{ij}_{st}\}_{s,t}$
 admits the representation via the potential
 (see (\ref{3.8})), namely, $g^{ij}_{st}=D^*G^{ij,st}$
  for all indices $i,j,s,t$ from 1 to  d.
  According  Lemma \ref{lem2},
  \[
 (g^{ij}_{st})_\e z_{ij}=( D^*G^{ij,st} )_\e z_{ij} 
  = D^*(\e^2 G^{ij,st}_\e z_{ij})-\e^2 G^{ij,st}_\e \cdot D z_{ij}-2\e(\div G^{ij,st})_\e \cdot \n z_{ij}
  \]
  \[
   = D^*(\e^2 G^{ij,st}_\e z_{ij})-2\e(\div G^{ij,st} \cdot e^k )_\e z_{ijk}
   -\e^2 G^{ij,st}_\e \cdot D z_{ij},
     \]
       where the each matrix
     $\{D^*(\e^2 G^{ij,st}_\e z_{ij})\}_{s,t}$ (without summing over $i,j$) has  the form of
          (\ref{3.15})  and therefore is solenoidal. Thereby,
     \[
 D^* R_\e^0=D^*(g^{ij}_\e z_{ij}) =  D_{st}\left( (g^{ij}_{st})_\e z_{ij}\right)=
-D_{st}\left(
2\e(\div G^{ij,st} \cdot e^k )_\e z_{ijk}
   +\e^2 G^{ij,st}_\e \cdot D z_{ij}\right), 
     \]
     where
     \[
    ( \div G^{ij,st} \cdot e^k )_\e z_{ijk}=(\pa_m G^{ij,st}_{km})_\e z_{ijk}=-(\pa_m G^{ij,km}_{st})_\e z_{ijk},
     \]
     \[
     G^{ij,st}_\e\cdot D z_{ij}=( G^{ij,st}_{km})_\e z_{ijkm}=-( G^{ij,km}_{st})_\e z_{ijkm},
     \]
thanks to the skew symmetry property of the matrix potential  $\{G^{ij,st}_{km}\}_{k,m}$.     Consequently,
 \begin{equation} \label{05.10}\ds{
 D^* R_\e^0=D_{st}\left(2\e (\pa_m G^{ij,km}_{st})_\e z_{ijk}+\e^2( G^{ij,km}_{st})_\e z_{ijkm}
 \right)}
 \atop\ds{
 =D^*
 \left(2\e (\pa_m G^{ij,km})_\e z_{ijk}+\e^2( G^{ij,km})_\e z_{ijkm}
 \right).}
\end{equation} 

From  (\ref{05.9}) and  (\ref{05.10}), we obtain the decomposition with only positive powers $\e^n$, $n\ge 1$: 
 \begin{equation} \label{05.11}
 \ds{
 D^* R_\e=\e D^*\left[
 (a(DN^{ijk}+2a(\n N^{ij}\times e^k)+2\pa_m G^{ij,km}_{st} - b^{ijk}
)_\e z_{ijk}\right]
} \atop\ds{
+\e^2 D^*\left[(aN^{ij})_\e Dz_{ij}+2 a_\e  (\n N^{ijk}\times e^m)_\e 
z_{ijkm}+( G^{ij,km})_\e z_{ijkm}\right]+
\e^3  D^*\left[(aN^{ijk})_\e Dz_{ijk}\right].
}
\end{equation}

In the sum (\ref{05.11}), the term with the first-order power $\e$ contains the oscillating function
\[
(a(DN^{ijk}+2a(\n N^{ij}\times e^k)+2\pa_m G^{ij,km}_{st} - b^{ijk}
)_\e\stackrel{(\ref{3.10})}=g^{ijk}_\e
\]
and can be written shortly as
 \begin{equation} \label{05.12} 
\e D^*[g^{ijk}_\e z_{ijk}]=\e D_{st}[(g^{ijk}_{st})_\e z_{ijk}].
\end{equation}
Here, the periodic matrix $g^{ijk}{=}\{g^{ij}_{st}\}_{s,t}
$
 possesses the properties  (\ref{ss}) and, thus, admits the representation, due to  Lemma \ref{lem1},
 via the matrix potential (see (\ref{3.13})$_1$), namely, $g^{ijk}_{st}=D^*G^{ijk,st}$
 for all indices $i,j,k,s,t$ from 1 to  d.
 By   Lemma \ref{lem2}, 
the  transformation is possible:
 \[
 (g^{ijk}_{st})_\e z_{ijk}=( D^*G^{ijk,st} )_\e z_{ijk} %
 \stackrel{(\ref{3.14})}
  = D^*(\e^2 G^{ijk,st}_\e z_{ijk})-\e^2 G^{ijk,st}_\e \cdot D z_{ijk}-2\e(\div G^{ijk,st})_\e \cdot \n z_{ijk}
  \]
  \[
   = D^*(\e^2 G^{ijk,st}_\e z_{ijk})-2\e(\div G^{ijk,st} \cdot e^p )_\e z_{ijkp}
   -\e^2 G^{ijk,st}_\e \cdot D z_{ijk}
     \]
    where the each matrix 
     $\{D^*(\e^2 G^{ijk,st}_\e z_{ijk})\}_{s,t}$ (without summing over $i,j,k$) has the form of
          (\ref{3.15})  and therefore is solenoidal.
 Thereby,  (\ref{05.12}) is rewritten as
\begin{equation} \label{05.13}
 \ds{
\e D^*[g^{ijk}_\e z_{ijk}]=2\e^2D_{st}
(\pa_q G^{ijk,pq}_{st})_\e z_{ijkp}
+\e^3D_{st}( G^{ijk,pq}_{st})_\e z_{ijkpq}
}\atop\ds{
=D^*[2\e^2 (\pa_q G^{ijk,pq})_\e z_{ijkp}+
\e^3 G^{ijk,pq}_\e z_{ijkpq}
],}
\end{equation}
where at the last step we simplify the write-up thanks to the skew symmetry property of the matrix potentials in the same way as in the case
of (\ref{05.10}).
   
   Eventually, the decomposition  (\ref{05.11}), with account of  (\ref{05.13}), runs as follows:
    \begin{equation} \label{05.14} 
  \ds{
 D^*R_\e =
 \e^2 D^*\left[2(\pa_q G^{ijk,pq})_\e z_{ijkp}+(aN^{ij})_\e Dz_{ij}+2 a_\e  (\n N^{ijk}\times e^m)_\e 
z_{ijkm}+( G^{ij,km})_\e z_{ijkm}\right]
}\atop\ds{
+\e^3  D^*\left[ G^{ijk,pq}_\e z_{ijkpq}+(aN^{ijk})_\e Dz_{ijk}\right]=D^*r_\e,
}
\end{equation}
where the term $r_\e$ coincides with (\ref{5.7}). Finally, (\ref{05.6}), (\ref{05.7}) together with (\ref{05.14}) entail the desired equality (\ref{5.5}).

To prove estimates (\ref{5.9}) and (\ref{5.90}), we apply
lemmas \ref{LemM1} and \ref{LemM5} while studying the terms $D^*r_\e$ and $r^0_\e$ in the right-hand side of (\ref{5.5}) (for more details see Section 6.2), and also the property of type (\ref{m.7}) for the smoothong operator $\Theta^\e $  while studying the term $(f^{,\e}-f)$.

From (\ref{5.5})--(\ref{5.7}) combined with the equality $f=(A_\e+1)u^\e$, we derive the equation solved by  $v^\e-u^\e$, namely, 
 \[
A_\e(\tilde{u}^\e-u^\e)+(\tilde{u}^\e-u^\e)=A_\e \tilde{u}^\e+\tilde{u}^\e-f
 \stackrel{(\ref{5.5})}=D^*r_\e+r^0_\e+(f^{,\e}-f)=F^\e.
\]
The energy estimate holds
\[
\|\tilde{u}^\e-u^\e\|_{ H^{2}(\rd)}\le c\|F^\e\|_{H^{-2}(\rd)},\quad c=const(d,\lambda),
\]
whereof    (\ref{2.17}) follows, due to (\ref{5.9}).
The proof of
Lemma \ref{Lem5.1} and Theorem \ref{th2.1} is complete.

\textbf{6.2.} Now, we clarify the step in the above proof connected with the derivation of
estimates (\ref{5.9}) and (\ref{5.90}), 
 technicalities of which have been formerly  dropped.
We observe, first, that $r_\e$ and $r^0_\e$ are the sums of six or two terms respectively (see (\ref{5.7}) and (\ref{5.6})). All these terms, except for the last two ones in (\ref{5.7}), can be estimated by Lemma \ref{LemM1}; as for the excluded ones from (\ref{5.7}), they have the common  factor $\e^3$, and Lemma \ref{LemM5} is employed to them. Here, the role of the oscillating multiplier $b(y)$ is played by one of the functions
 \begin{equation} \label{05.18}
N^{ij},\, N^{ijk},\, \n N^{ijk},\, G^{ij,kp},\, G^{ijk,pq},\,\n G^{ijk,pq}
\end{equation}
 with indices $i,j,k,p, q$ from 1 to d. All these periodic functions belong to the space $L^2_\per(Y)$,
which stems from the energy estimate (\ref{3.4}) for $N^{ij}$; moreover, by any choice of the  multiplier $b$
from the set of functions (\ref{05.18}),  the inequality $\langle|b|^2\rangle\le C$ is valid with constant $C$ depending only on
 $d$ and $\lambda$.

Recalling the notation (\ref{5.8}) for the functions   
 $z_{ij}$, $z_{ijk}$, $z_{ijkp}$ participating  in (\ref{5.7}) and (\ref{5.6}), we see that each of them acquires the form 
 $\Theta^\e \varphi$ with 
$\varphi\in \ld$ satisfying 
the inequality
$\|\varphi\|_{\ld}\le \|\hat{u}^\e\|_{H^4(\rd)}\le C$ 
 in view of (\ref{2.13}). Therefore, the estimate         of the type (\ref{m.4}) with the smoothing operator $\Theta^\e$ can be applied to the terms containing
              $z_{ij}$, $z_{ijk}$, $z_{ijkp}$. For instance, since  $z_{ijkp}=\Theta^\e D_iD_jD_kD_p
\hat{u}^\e$, 
it follows that
\[
\|( G^{ij,kp})_\e z_{ijkp}\|\stackrel{(\ref{m.4})}\le c\langle |G^{ij,kp}|^2\rangle^{1/2}
\| D_iD_jD_kD_p \hat{u}^\e\|
\stackrel{(\ref{2.12})}\le 
C\|f\|,
\]
where at the last step                                                                                                                                                                                                                                                                                                                                                                                                                                                                                                                                                                                                                                                                                                                             the  $L^2$-estimate for the periodic function $G^{ij,kp}$ is exploited.

The terms from (\ref{5.7})  containing $z_{ijkpq}$ or $D z_{ijk}$ can be estimated by Lemma \ref{LemM5} with sacrifice of  the factor $\e$ 
as a payment for the  uniform   with respect to $\e$  bound.
As  a result, 
these terms, though including the 
 power $\e^3$,  really are of the order
  $O(\e^2)$.
For instance,  $z_{ijkpq}=\Theta^\e(D_q \varphi)$, where
$\varphi= D_iD_jD_kD_p \hat{u}^\e$, thereby,
 \begin{equation} \label{6.10}
\e \|( G^{ijk,pq})_\e z_{ijkpq}\|\stackrel{(\ref{m.13})}\le C\langle |G^{ijk,pq}|^2\rangle^{1/2}
\| D_iD_jD_kD_p\hat{u}^\e\|
\stackrel{(\ref{2.12})}\le 
C\|f\|,
\end{equation}
because  the smoothing kernel of the operator $\Theta^\e$   is sufficiently regular for Lemma \ref{LemM5} to be applied here. At the last step                                                                                                                                                                                                                                                                                                                                                                                                                                                                                                                                                                                                                                                                                                                            
in (\ref{6.10}), the  $L^2$-estimate for the periodic function $G^{ijk,pq}$ is exploited too.

\section{$\!\!\!\!\!\!$. Some remarks}\label{sec7}
 \setcounter{theorem}{0} \setcounter{equation}{0}
  $\quad{}$ \textbf{Remark 1.}  Our proof of estimates (\ref{2.17}) and (\ref{2.22}) does not rely in an essential way on the selfadjointness of the operator $A_\e$. Similar estimates hold true if the symmetry condition from the assumption (\ref{2.2}) is omitted.
The advantage of  the symmetry condition is that our formulas become less cumbersome, and a lower number of objects is  requested for their derivation. 
   For instance,  passing from (\ref{5.3}) to (\ref{5.4}), we rewrite 
   the sum of five terms as the sum of only four terms. Moreover, the symmetry property  of the tensor $a(y)$ is inherited by other objects of our consideration, e.g., 
$g^{ij}$, $G^{ij,st}$, $g^{ijk}$, $G^{ijk,st}$ (see \S3);  as a consequence,  the formulas, where these 
 objects are involved,  may  also become simpler. In particular,  the equality  (\ref{3.16}) (employed not once in the representation of the type (\ref{3.14})) has a more cumbersome version in the case of a nonsymmetric matrix $B$. 

\textbf{Remark 2.} Assume that the symmetry condition is not required in (\ref{2.2}). Then the $H^2$-approximation of $\e^2$ order for the resolvent of  $A_\e$ constructed in a similar fashion as in (\ref{2.13})--(\ref{2.16}) provides at the same time the 
$L^2$-approximation of $\e^2$ order and, thus, 
$$(A_\e+1)^{-1}=(\hat{A}_\e+1)^{-1}+O(\e^2),$$
if the correctors are dropped and the main term is simplified by properties of smoothing. The latter approximation has a right to be mentioned  along with other ones, e.g., that of  (\ref{07}).

\textbf{Remark 3.}  The approximations (\ref{08}) and (\ref{2.21}) of the resolvent $A_\e$ involve
 the resolvent of the homogenized fifth-order operator 
   $\hat{A}_\e$,  although the order of $A_\e$ is four. This phenomenon is known;   in various homogenization problems, there appeared formerly not once homogenized operators, which had the order larger than that of the given operator.
For instance, 
  in \cite{P12}, where a   non-stationary diffusion equation in a periodic medium is considered
  and approximations for its solution  $u(x,t)$ are sought in
 $L^2$-norm in  sections  $t{=}const{>}0$ with accuracy $O(t^{-\f{m}{2}})$ 
($m{\ge} 2$ is a natural parameter) as $t{\to} +\infty$, 
the order of the  homogenized operator is dictated by the order of approximation and increases without limit simultaneously with the parameter $m$.

 \textbf{Remark  4.} The operator $\hat{A}_\e$ appears as a singular perturbation of the   operator $\hat{A}$ commonly employed in homogenization.
 Singularly perturbed operators 
 are known and widely used in the elasticity theory of thin plates and shells
 (see, e.g., \cite{G}). But the  
 operator $\hat{A}_\e$, which is involved in 
 (\ref{2.13}) or (\ref{2.21}), does not  refer to any type of the singular perturbation described in
  \cite{G}.
 
  \textbf{Remark 5.}
 In the proof of Theorem \ref{th5.1}, any expression  $T$ admitting the estimate
 $|T|\le c\e^3\|f\|\,\|h\|$ is considered to be nonessential. Some terms classified as nonessential in this sense actually satisfy sharper estimates, e.g., with  majorants of $\e^4$ order. For example, in Section  5.3 at the step 2$^\circ$, all the  terms declared as nonessential are in fact of this kind.
Consider one of them:
\[
\e^2 (N^{ij}_\e z_{ij}, v^{,\e})\stackrel{(\ref{5.8}), (\ref{5.17})}=\e^2 (N^{ij}_\e S^\e \varphi, S^\e\psi)
\stackrel{(\ref{m.9})}\le C\e^4\langle |N^{ij}|^2\rangle^{1/2} \|\n\varphi\|\,\| \nab\psi\|,
\]
by Lemma \ref{LemM3}, which can be applied here because $\varphi=D_iD_j(S^\e)^2\hat{u}^\e$ and $\psi=(S^\e)^2\hat{v}^\e$ belong to $H^1(\rd)$ with an appropriate estimate of the norms. Such kind of refinement in the estimates 
is necessary if   $\e^4$ order approximations are sought.

\textbf{Remark 6.} Some words about the method we use to prove operator estimates. We apply 
 the approach proposed firstly by Zhikov in \cite{Zh1},  
but in the version from \cite{ZhP05}. While constructing approximations for the resolvent of the operator $A_\e$ and then justifying them, we encounter serious problems linked with the minimal regularity in data when even the correctness of these approximations is under the question. To cope with this kind of difficulties, there are 
two close to each other ways proposed in \cite{Zh1} and  \cite{ZhP05}. The first way is to introduce a pure shift parameter into the problem (more exactly, into the coefficients of the operator,  and this perturbation is inherited in  the approximations)
with further integrating over the  shift parameter
(as, e.g., in \cite{Zh1} or \cite{PMSb}); the second way is to insert smoothing operators inside the approximations themselves from the very beginning (as in \cite{ZhP05} or \cite{P20a}--\cite{P20s}).

Another method  to prove operator estimates in homogenization is based on the Floquet--Bloch transformation and spectral arguments. It is used for high order operators, e.g., in  \cite{V},  \cite{KS},  and \cite{SS}.
   The result recently announced in  \cite{SS} is related to a fourth-order elliptic self-adjoint matrix differential operator $A_\varepsilon$ with $\varepsilon$-periodic coefficients which satisfies a certain factorization condition. 
 In   \cite{SS}, we find resolvent 
   approximations of the same order of smallness in error as in Theorem \ref{th5.1}, but they look rather differently compared to that of (\ref{2.22}).


\end{document}